\documentclass[10pt,a4paper]{article}
\usepackage{amsmath,amssymb,amsthm,graphics,graphicx,epic,eepic,multicol,ascmac}
\usepackage[mathscr]{euscript}
\usepackage[all]{xy}
\usepackage{url}

\numberwithin{equation}{section}
\theoremstyle{theorem}
\newtheorem{thm}{Theorem}[section]
\newtheorem{prop}[thm]{Proposition}
\newtheorem{lem}[thm]{Lemma}
\newtheorem{rem}[thm]{Remark}
\newtheorem{cor}[thm]{Corollary}

\theoremstyle{definition}
\newtheorem{defn}[thm]{Definition}
\newtheorem{ex}[thm]{Example}

\def\al{\alpha}

\def\wht(#1){\widehat{\ #1\ }}

\newcommand{\cA}{{\mathcal A}}

\newcommand{\cX}{{\mathcal X}}

\newcommand{\frg}{\mathfrak g}
\newcommand{\frh}{\mathfrak h}

\newcommand{\frl}{\mathfrak l}

\newcommand{\frs}{\mathfrak s}

\newcommand{\bbC}{\mathbb C}

\newcommand{\bbZ}{\mathbb Z}

\newcommand{\Ad}{\mathrm{Ad}}

\newcommand{\ch}{\mathrm{ch}}

\newcommand{\lbr}{\begin{bmatrix}}
\newcommand{\rbr}{\end{bmatrix}}

\newcommand{\cd}{commutative diagram }

\def\ge{\frg}

\def\al{\alpha}

\def\beneme{\begin{enumerate}}
\def\beq{\begin{equation}}
\def\beqn{\begin{eqnarray}}
\def\beqnn{\begin{eqnarray*}}

\def\bfii0{{\bf i_0}}

\def\bbra#1,#2,#3{\left\{\begin{array}{c}\hspace{-5pt}
#1;#2\\ \hspace{-5pt}#3\end{array}\hspace{-5pt}\right\}}
\def\cd{\cdots}
\def\ci(#1,#2){c_{#1}^{(#2)}}
\def\Ci(#1,#2){C_{#1}^{(#2)}}
\def\mpp(#1,#2,#3){#1^{(#2)}_{#3}}
\def\bCi(#1,#2){\ovl C_{#1}^{(#2)}}
\def\ch(#1,#2){c_{#2,#1}^{-h_{#1}}}
\def\cc(#1,#2){c_{#2,#1}}

\def\Del{\Delta}

\def\di(#1,#2){D_{#1}^{(#2)}}
\def\dbi(#1,#2){\ovl D_{#1}^{(#2)}}

\def\eneme{\end{enumerate}}

\def\eeq{\end{equation}}
\def\eeqn{\end{eqnarray}}
\def\eeqnn{\end{eqnarray*}}

\def\gau#1,#2{\left[\begin{array}{c}\hspace{-5pt}#1\\
\hspace{-5pt}#2\end{array}\hspace{-5pt}\right]}

\def\ji(#1,#2){j_{#1}^{(#2)}}

\def\lan{\langle}

\def\lm{\lambda}
\def\Lm{\Lambda}

\def\nd{\noindent}

\def\ovl{\overline}

\def\qq{\qquad}
\def\q{\quad}
\def\qed{\hfill\framebox[2mm]{}}

\def\ran{\rangle}

\def\TY(#1,#2,#3){#1^{(#2)}_{#3}}

\def\xxi(#1,#2,#3){\displaystyle {}^{#1}\Xi^{(#2)}_{#3}}
\def\xsi(#1,#2,#3){\displaystyle {}^{#1}\Sigma^{(#2)}_{#3}}
\def\xE(#1,#2,#3){\displaystyle {}^{#1}E_{#2}[#3]}
\def\xF(#1,#2){\displaystyle {}^{#1}F_{#2}}
\def\xx(#1,#2){\displaystyle {}^{#1}\Xi_{#2}}
\def\W1{W(\varpi_1)}

\def\m@th{\mathsurround=0pt}
\def\fsquare(#1,#2){
\hbox{\vrule$\hskip-0.4pt\vcenter to #1{\normalbaselines\m@th
\hrule\vfil\hbox to #1{\hfill$\scriptstyle #2$\hfill}\vfil\hrule}$\hskip-0.4pt
\vrule}}

\newcommand{\ba}{\begin{array}}
\newcommand{\ea}{\end{array}}

\newcommand{\eq}{\begin{eqnarray}}
\newcommand{\eneq}{\end{eqnarray}}

\title{\textbf{\large{Geometric crystals and Cluster ensembles in Kac-Moody setting}}}
\author{\normalsize{YUKI KANAKUBO\thanks{Division of Mathematics, 
Sophia University, Kioicho 7-1, Chiyoda-ku, Tokyo 102-8554,
Japan: {j\_chi\_sen\_you\_ky@eagle.sophia.ac.jp. }}
\ and\ 
TOSHIKI NAKASHIMA\thanks{Division of Mathematics, 
Sophia University, Kioicho 7-1, Chiyoda-ku, Tokyo 102-8554,
Japan: {toshiki@sophia.ac.jp}.}
}}
\date{}

\begin{document}

\maketitle
\vspace{-10pt}

\begin{abstract}
For a Kac-Moody group $G$, double Bruhat cells $G^{u,e}$ ($u$ is a Weyl group element)
have positive geometric crystal structures. 
In \cite{HW}, it is shown that there exist birational maps 
between `cluster tori' $\mathcal{X}_{\Sigma}$ (resp. $\mathcal{A}_{\Sigma}$) and
$G_{\Ad}^{u,e}$ (resp. $G^{u,e}$), and they are extended
to regular maps from cluster $\mathcal{X}$ (resp. $\mathcal{A}$)
-varieties to $G_{\Ad}^{u,e}$ (resp. $G^{u,e}$).
The aim of this article is to construct certain positive 
geometric crystal structures on the cluster tori $\cX_{\Sigma}$ and $\cA_{\Sigma}$
by presenting their explicit formulae. In particular,
the geometric crystal structures on the tori $\cA_{\Sigma}$ are obtained
by applying the twist map. 
As a corollary, we see the sets of $\mathbb{Z}^T$-valued points of the cluster varieties
have plural structures of crystals.
\end{abstract}

\section{Introduction}

In \cite{FG},
V.V.Fock and A.B.Goncharov have introduced the ``cluster ensemble'', 
which is a 
pair of schemes ``$\cA$-variety $\cA_{|\Sigma|}$" and ``$\cX$-variety $\cX_{|\Sigma|}$''
with a morphism $\cA_{|\Sigma|}\rightarrow \cX_{|\Sigma|}$ called an {\it ensemble map}. 
For the construction of these schemes, {\it cluster} $\mathcal{A}$ (resp. $\mathcal{X}$)-{\it tori}
$\cA_\Sigma=\{(A_i)\,|\, i\in I,\ A_i\in\bbC^\times\}$
(resp. $\cX_\Sigma=\{(X_i)\,|\, i\in I,\ X_i\in\bbC^\times\}$)
associated with {\it seeds} $\Sigma=(I,I_0,B,d)$ (see Sect.\ref{ensect}) are defined. 
The schemes are defined as 
unions of cluster tori which are glued by the following 
{\it mutations} $\mu_k:\cA_{\Sigma}\rightarrow \cA_{\Sigma'}$, $\cX_{\Sigma}\rightarrow \cX_{\Sigma'}$
\begin{eqnarray*}
&&\mu_k^*(A'_{M_k(i)})=
\begin{cases}
A_i & {\rm if}\ i\neq k,\\
\frac{\prod_{b_{k,j}>0}A_j^{b_{k,j}}+\prod_{b_{k,j}<0}A_j^{-b_{k,j}}}{A_k} & {\rm if}\ i=k,
\end{cases}
\\
&&\mu_k^*(X'_{M_k(i)})=
\begin{cases}
X_iX_k^{[b_{i,k}]_+}(1+X_k)^{-b_{i,k}} & {\rm if}\ i\neq k,\\
X_k^{-1} & {\rm if}\ i=k,
\end{cases}
\end{eqnarray*}
where $\Sigma'=(I',I_0,B',d')$ is a new seed, and the map $M_k:I\rightarrow I'$ 
is as in Sect.\ref{ensect}.


Furthermore, in \cite{FG} 
Fock and Goncharov presented a conjecture on ``tropical duality" between 
these two cluster varieties. To be more precise, the conjecture claimed that 
the ``universal positive Laurent polynomial ring'' on 
$\cA_{|\Sigma|}$(resp. $\cX_{|\Sigma|}$) is described by the positive summation of 
points in the set of $\mathbb{Z}$-valued points (see \ref{zvalue}) $\cX_{|\Sigma^{\vee}|}(\bbZ^T)$
(resp. $\cA_{|\Sigma^{\vee}|}(\bbZ^T)$), where $\Sigma^{\vee}$ is the Langlands dual seed
of $\Sigma$ \cite{FG}.
Though in \cite{GHK}, unfortunately, the counter-examples for the conjecture 
have been found, and in \cite{GHKK}, the conjecture has been refined 
to be more valid, which is now called ``full Fock-Goncharov conjecture".
It seems to be still generically open except for several special cases.

H.Williams considered the relations between cluster ensembles associated with 
reduced words of Weyl group elements and 
double Bruhat cells associated with those Weyl group elements 
in Kac-Moody setting and furthermore, he also gave the isomorphism 
between the coordinate ring $\bbC[\cA_{|\Sigma|}]$ and that on the double Bruhat 
cell (\cite{HW}). 
Therein, he also 
constructed the regular map $p_M$ from $\cA_{|\Sigma|}$ to $\cX_{|\Sigma|}$
compatible with all cluster mutations (Proposition \ref{map-p}),
which plays a role of the ensemble map in the context of \cite{HW}. 

In \cite{BK}, A.Berenstein and D.Kazhdan 
have initiated the theory of ``geometric crystal", which is aimed 
to construct a geometric analogue of the Kashiwara's crystal base theory
on a variety birationally isomorphic to a split torus. 
In \cite{N0}, the second author extended this notion to the Kac-Moody setting
and gave some explicit forms of geometric crystals on Schubert/Bruhat cells.
Geometric crystals have a bunch of remarkable properties, in particular, 
the fact that positive geometric crystals can be transferred to the Langlands 
dual Kashiwara's crystal bases by the ``tropicalization" procedure 
is one of the most crucial features. More precisely, there exists
a functor $\mathcal{U}\mathcal{D}$ from a category of split tori to the category of set,
and each variety which has a positive geometric crystal structure
 corresponds to the set of its co-characters under $\mathcal{U}\mathcal{D}$, which
has a crystal structure.

The aim of this article is to construct certain positive 
geometric crystal structures on the cluster tori $\cX_{\Sigma}$ and $\cA_{\Sigma}$ of \cite{HW}
by presenting their explicit formulae.
As a corollary of them, we see that
the sets of $\mathbb{Z}^T$-valued points $\cX_{|\Sigma|}(\bbZ^T)$, $\cA_{|\Sigma|}(\bbZ^T)$
have {\it glued crystal} structures.
Here, glued crystal is defined as a unit of some crystals (see \ref{glued-cry-sec}).
We expect those crystals will be a guide to resolve the Fock-Goncharov 
conjectures. 
To achieve the aim, first we define geometric crystal structures on the
cluster tori $\cX_{\Sigma}$, $\cA_{\Sigma}$. 
We will see that double Bruhat cells $G^{u,e}$ and their quotients
 $G^{u,e}_{\rm Ad}$ have geometric crystal structures in Proposition \ref{GCcell} and Definition \ref{GCdefX}.
Using birational maps from $\cA_{\Sigma}$, $\cX_{\Sigma}$ to $G^{u,e}$, $G^{u,e}_{\rm Ad}$
given by H.Williams, we obtain geometric crystal structures on
$\cA_{\Sigma}$, $\cX_{\Sigma}$ (Definition \ref{actions}(1)).
By using ``twist map" $\zeta^{u,e}$, we can also construct 
another geometric crystal structures on them (Definition \ref{actions}(2)).
Second, we will verify compatibilities between these structures
in Proposition \ref{compprop}.
We mainly treat the geometric crystal structures of Definition \ref{actions}(1) on $\cX_{\Sigma}$ and
those of Definition \ref{actions}(2) on $\cA_{\Sigma}$ in this article. 
Third, we will present explicit formulae of geometric crystal
structures on the tori $\cA_{\Sigma}$, $\cX_{\Sigma}$.
These explicit formulae imply the tori
have positive geometric crystal structures.
Since the set of $\mathbb{Z}^T$-valued points $\cA_{|\Sigma|}(\bbZ^T)$
 (resp. $\cX_{|\Sigma|}(\bbZ^T)$) is a union of the sets of co-characters of $\cA_{\Sigma}$ (resp. $\cX_{\Sigma}$),
we see that it has a glued crystal structure.

The organization of this article is as follows. In Sect.2, we review on 
the theory of geometric crystals and some explicit formula of 
geometric crystal on some double Bruhat cell $G^{u,e}$ which 
will be needed in the rest of the article.
In Sect.3, the cluster ensembles will be introduced and we will see the
cluster ensembles associated with arbitrary Weyl group elements. 
A cluster tori $\cA_{\Sigma_{\rm \bf {i}}}$, $\cX_{\Sigma_{\rm \bf {i}}}$
birationally isomorphic to $G^{u,e}$, $G^{u,e}_{\rm Ad}$
are defined in this section.
In Sect.4, we define geometric crystal structures
on the cluster tori $\cA_{\Sigma_{\rm \bf {i}}}$, $\cX_{\Sigma_{\rm \bf {i}}}$.
We will also present a compatibility between
these structures.
In Sect.5, we give explicit formulae of
geometric crystal structures on $\cX_{\Sigma_{\rm \bf {i}}}$.
In Sect.6, we also give explicit formulae of
geometric crystal structures on $\cA_{\Sigma_{\rm \bf {i}}}$
by using the compatibility shown in Sect.4.
In Sect.7, we present explicit formulae of
geometric crystal structures on $\cA_{\Sigma_{\rm \bf {i}}}$
in the case $G=SL_{r+1}(\mathbb{C})$ and $u$ is the longest element of $W$
in a different way from Sect.6.

{\bf Acknowledgements} We would like to appreciate David Hernandez, 
Michael Gekhtman, Gleb Koshevoy, 
Beranrd Leclerc, and Hironori Oya
for fruitful discussions and valuable suggestions. 
We also thank the faculty and staffs at Universit'e Paris Diderot during
their stay in 2017-2018.
and T.N. thanks the faculty and staffs at IHES for his stay in 2017. 
Y.K. was supported by JSPS KAKENHI Grant Number 17H07103, 
and T.N. was supported in part by JSPS Grants in Aid for Scientific Research $\sharp$15K04794.

\section{Geometric crystals}

In this section, we will review notion of geometric crystal following \cite{BK0, BK, N0,N}.

\subsection{Notation and definitions}\label{notation}

Following \cite{HW}, we define several notion.
For a positive integer $l$, we set $[1,l]:=\{1,2,\cdots,l\}$.
Let $A=(a_{i,j})_{i,j\in [1,r]}$ be a symmetrizable generalized
Cartan matrix with a symmetrizer ${\rm diag}(d_1,\cdots,d_r)$ $(d_i\in\mathbb{Z}_{>0})$,
and $\frg=\frg(A)=\lan e_i,f_i,\frh \ran$ the Kac-Moody Lie algebra associated with $A$ over $\mathbb{C}$.
The Cartan subalgebra $\frh\subset\frg$ contains simple coroots
$\al_1^{\vee},\cdots,\al_r^{\vee}$ and $\frh^*$ contains simple roots $\al_1,\cdots,\al_r$,
which satisfy $\al_j(\al_i^{\vee})=a_{i,j}$. Let $\tilde{r}={\rm dim}\ \frh=2r-{\rm rank}\ A$.
The simple reflections $s_i\in {\rm Aut}(\frh^*)$ $(i\in [1,r])$ are defined as
$s_i(\beta):=\beta-\beta(\al_i^{\vee})\al_i$, which generate the Weyl group $W$.
Let $P:=\{\lambda\in\frh^* | \lambda(\al_i^{\vee})\in\mathbb{Z}\ {\rm for\ all}\ i\in [1,r]\}$ and
$\{\Lm_i\}_{i \in [1,\tilde{r}]}\subset P$ be a basis which satisfies
$\Lm_i(\al_j^{\vee})=\delta_{i,j}$ for $i\in[1,\tilde{r}]$, $j\in[1,r]$.
Let $G$ be the Kac-Moody group associated with $(\frg,P)$ and $H\subset G$ a maximal torus.
The set $P$ is identified with the set
${\rm Hom}(H,\mathbb{C}^{\times})$ of characters.
We call $\Lm_i$ $(i \in [1,\tilde{r}])$ {\it fundamental weights}.
By fixing $\{\Lm_i\}_{i \in [1,\tilde{r}]}$,
we obtain a corresponding dual basis of ${\rm Hom}(\mathbb{C}^{\times},H)$ and denote its
elements $\al^{\vee}_1,\cdots,\al^{\vee}_{\tilde{r}}$.
If $i\in[1,r]$ then $\al^{\vee}_i$ is just the $i$-th coroot of $\frg$.
For $t\in\mathbb{C}^{\times}$ and $h\in{\rm Hom}(\mathbb{C}^{\times},H)$,
let $t^{h}$ denote the element $h(t)\in H$.
The equation
\begin{equation}\label{al-rel}
\al_j=\sum_{1\leq i\leq \tilde{r}} a_{ij}\Lm_i
\end{equation}
defines numbers $a_{ij}$ for $i\in\{r+1,\cdots,\tilde{r}\}$ and $j\in\{1,2,\cdots,r\}$.
We can also define elements $\{\al_i\}_{i\in\{r+1,\cdots,\tilde{r}\}}$ of $P$ by
\[
\al_i=D\sum^{r}_{j=1} d_j^{-1} a_{ij}\Lm_j,
\]
where $D$ is the least common integer multiple of 
$d_1,\cdots,d_r$. 

For each real root $\al$,
there exists a one-parameter subgroup $\{x_{\al}(t)|t\in\mathbb{C}\}\subset G$,
and $G$ is generated by all one-parameter subgroups and $H$ \cite{Kum,KP}.
Let $N$, $N^-$ be the subgroups of $G$
generated by $\{x_{\al}(t)|t\in\mathbb{C},\ \al:{\rm positive\ root}\}$, $\{x_{\al}(t)|t\in\mathbb{C},\ \al:{\rm negative\ root}\}$.
Let $B=HN$, $B^-=HN^-$ be Borel subgroups.
For $T\in H$, let $\al_i(T)$ or $T^{\al_i}$ denote the value of $\al_i$ at $T$ (as the character).

\subsection{Double Bruhat cells}

We set $x_i(c):={\rm exp}(c e_i)$, $y_i(c):={\rm exp}(c f_i)\in G$ for $c\in\mathbb{C}$.
We also set $\ovl{s_i}:=x_{i}(-1)y_i(1)x_i(-1)$ for $i\in[1,r]$, and for a reduced expression $w=s_{j_1}\cdots s_{j_n}\in W$, set $\ovl{w}:=\ovl{s_{j_1}}\cdots \ovl{s_{j_n}}$. The following two kinds of Bruhat decompositions of $G$ are known \cite{Kum}:
\[ G=\displaystyle\coprod_{u \in W}B\ovl{u}B=\displaystyle\coprod_{u \in W}B^-\ovl{u}B^- .\]
Then, for $u$, $v\in W$, 
the {\it double Bruhat cell} $G^{u,v}$ is defined as follows:
\[ G^{u,v}:=B\ovl{u}B \cap B^-\ovl{v}B^-. \]

\begin{prop}{\rm \cite{HW}}
For $u$, $v\in W$, the double Bruhat cell $G^{u,v}$ is a rational affine variety and
{\rm dim} $G^{u,v}=l(u)+l(v)+\tilde{r}$.
\end{prop}

\subsection{Crystals}

Let us recall the definition of {\it crystals} \cite{K2}. 
We use the notation in \ref{notation}.

\begin{defn}
A {\it crystal} is a set $\mathcal{B}$ together with the maps
${\rm wt}_i:\mathcal{B}\rightarrow \mathbb{Z}$,
$\varepsilon_i,\varphi_i:\mathcal{B}\rightarrow \mathbb{Z}\cup \{-\infty\}$
and $\tilde{e}_i$,$\tilde{f}_i:\mathcal{B}\rightarrow \mathcal{B}\cup\{0\}$
($i\in[1,r]$) satisfying the followings: For $b,b'\in\mathcal{B}$, $i,j\in[1,r]$,
\begin{enumerate}
\item[$(1)$] $\varphi_i(b)=\varepsilon_i(b)+{\rm wt}_i(b)$,
\item[$(2)$] ${\rm wt}_j(\tilde{e}_ib)={\rm wt}_j(b)+a_{j,i}$ if $\tilde{e}_i(b)\in\mathcal{B}$,
\quad ${\rm wt}_j(\tilde{f}_ib)={\rm wt}_j(b)-a_{j,i}$ if $\tilde{f}_i(b)\in\mathcal{B}$,
\item[$(3)$] $\varepsilon_i(\tilde{e}_i(b))=\varepsilon_i(b)-1,\ \ 
\varphi_i(\tilde{e}_i(b))=\varphi_i(b)+1$\ if $\tilde{e}_i(b)\in\mathcal{B}$, 
\item[$(4)$] $\varepsilon_i(\tilde{f}_i(b))=\varepsilon_i(b)+1,\ \ 
\varphi_i(\tilde{f}_i(b))=\varphi_i(b)-1$\ if $\tilde{f}_i(b)\in\mathcal{B}$, 
\item[$(5)$] $\tilde{f}_i(b)=b'$ if and only if $b=\tilde{e}_i(b')$,
\item[$(6)$] if $\varphi_i(b)=-\infty$ then $\tilde{e}_i(b)=\tilde{f}_i(b)=0$.
\end{enumerate}
We call $\tilde{e}_i$,$\tilde{f}_i$ {\it Kashiwara operators}, and ${\rm wt}_i$ {\it weight functions}.
A crystal $\mathcal{B}$ is said to be {\it free} if the Kashiwara operators
$\tilde{e}_i$ $(i\in [1,r])$ are bijections $\tilde{e}_i:\mathcal{B}\rightarrow \mathcal{B}$. 
\end{defn}

Note that the above definition of crystals is slightly weaker than the original one in \cite{K2}.
In the case the generalized Cartan matrix $(a_{i,j})_{i,j\in[1,r]}$ has rank $r$
then the above definition is equivalent to the one in \cite{K2}.

For two crystals $(\mathcal{B},\{\tilde{e}_i\},\{\tilde{f}_i\},\{\varepsilon_i\},\{\varphi_i\},\{{\rm wt}_i\})$,
$(\mathcal{B}',\{\tilde{e}_i'\},\{\tilde{f}_i'\},\{\varepsilon'_i\},\{\varphi'_i\},\{{\rm wt'}_i\})$,
a bijection $f:\mathcal{B}\rightarrow \mathcal{B}'$
is called a {\it crystal isomorphism} if it satisfies
$f(\tilde{e}_i(b))=\tilde{e}_i'(f(b))$, 
$f(\tilde{f}_i(b))=\tilde{f}_i'(f(b))$, 
$\varepsilon_i'(f (b))=\varepsilon_i(b)$,
$\varphi_i'(f(b))=\varphi_i(b)$ and ${\rm wt}_i'(f(b))={\rm wt}_i(b)$ for $b\in\mathcal{B}$ and $i\in [1,r]$.
Here we understand $f(0)=0$.

\subsection{Geometric crystals}\label{GCsub}

For algebraic varieties $X$, $Y$ and a rational function $f:X\rightarrow Y$,
let ${\rm dom}(f)$ denote the maximal open subset of
$X$ on which $f$ is defined.

\begin{defn}
For a symmetrizable generalized
Cartan matrix $A=(a_{i,j})_{i,j\in [1,r]}$ and
an irreducible algebraic variety $X$ over $\mathbb{C}$, let $\gamma_i$, $\varepsilon_i$ $(i\in [1,r])$ be rational functions on $X$, and $e_i : \mathbb{C}^{\times}\times X\rightarrow X$ a rational $\mathbb{C}^{\times}$-action $(i\in [1,r])$
(to be denoted by $(c,x)\mapsto e^c_i(x)$). A quintuple $(X,\{e_i\}_{i\in [1,r]},\{\gamma_i\}_{i\in [1,r]},\{\varepsilon_i\}_{i\in [1,r]})$ is called a {\it geometric crystal} if 
\begin{enumerate}
\item[$(i)$] For $i\in [1,r]$, $(\{1\}\times X)\cap {\rm dom}(e_i)$ is open dense in $\{1\}\times X$.
\item[$(ii)$] For any $i$, $j\in [1,r]$, the rational functions $\{\gamma_i\}_{i\in [1,r]}$ satisfy $\gamma_j(e^c_i(x))=c^{a_{ij}}\gamma_j(x)$.
\item[$(iii)$] For any $t\in H$, $w\in W$ and its two reduced words $\textbf{i}$, $\textbf{i}'$, the relation $e_{\textbf{i}}(t)=e_{\textbf{i}'}(t)$ holds, where for a reduced word $\textbf{i}=(i_1,\cdots,i_n)$ of $w$, we define $e_{\textbf{i}}(t)=e^{(\al^{(1)}(t))}_{i_1}e^{(\al^{(2)}(t))}_{i_2}\cdots e^{(\al^{(n)}(t))}_{i_n}$,\ \ $\al^{(j)}:=s_{i_n}\cdots s_{i_{j+1}}(\al_{i_j})$ .
\item[$(iv)$] The rational functions $\{\varepsilon_i\}_{i\in [1,r]}$ satisfy $\varepsilon_i(e^c_i(x))=c^{-1}\varepsilon_i(x)$ and $\varepsilon_i(e^c_j(x))= \varepsilon_i(x)$ if $a_{i,j}=a_{j,i}=0$.
\end{enumerate}
\end{defn}


Let $X^{*}(T):={\rm Hom}(T,\mathbb{C}^{\times})$
be the set of characters for
a split algebraic torus $T$.

\begin{defn}
Let $T$, $T'$ be split algebraic tori over $\mathbb{C}$.
\begin{enumerate}
\item[(i)] A regular function $f =\sum_{\mu\in X^{*}(T)} c_{\mu}\cdot \mu$ on $T$ is positive if all coefficients $c_{\mu}$ are non-negative
numbers. A rational function on $T$ is said to be positive if there exist positive regular functions
$g$, $h$ such that $f=\frac{g}{h}$ $(h\neq0)$.
\item[(ii)] Let $f : T \rightarrow T'$ be a rational map between $T$ and $T'$. Then $f$ is called {\it positive} if for any
$\xi\in X^*(T')$, the rational function $\xi\circ f$ is positive in the sense of (i).
\end{enumerate}
\end{defn}

Let $\mathcal{T}_+$ be a category whose objects are algebraic tori over $\mathbb{C}$ and morphisms are positive rational maps.
In \cite{BK0, N0}, a functor $\mathcal{U}\mathcal{D}:\mathcal{T}_+\rightarrow \mathfrak{S}\mathfrak{e}\mathfrak{t}$ is introduced, where $\mathfrak{S}\mathfrak{e}\mathfrak{t}$ is the category of all sets. Each torus $T$ in $\mathcal{T}_+$ corresponds to the set of
co-characters $X_{*}(T)={\rm Hom}(\mathbb{C}^{\times},T)$ under the functor $\mathcal{U}\mathcal{D}$.

\begin{defn}\label{pos-str}
Let $\chi=(X,\{e_i\}_{i\in [1,r]},\{\gamma_i\}_{i\in [1,r]},\{\varepsilon_i\}_{i\in [1,r]})$
be a geometric crystal, $T$ an algebraic torus, $\theta:T\rightarrow X$ a birational map.
The map $\theta$ is called {\it positive structure} on $\chi$ if it satisfies the following:
\begin{enumerate}
\item[(i)] For $i\in[1,r]$, the rational functions $\gamma_i\circ\theta$, $\varepsilon_i\circ\theta$ are positive.
\item[(ii)]  For $i\in[1,r]$, the rational map $e_{i,\theta}:\mathbb{C}^{\times}\times T\rightarrow T$,
defined by $(c,t)\mapsto \theta^{-1}\circ e_i^c\circ \theta(t)$ is positive.
\end{enumerate}
We say $(\chi,\theta)$ is a {\it positive geometric crystal}. 
\end{defn}

Applying the functor $\mathcal{U}\mathcal{D}$ to $e_{i,\theta}$, $\gamma_i\circ\theta$ and $\varepsilon_i\circ\theta$, we get
\[
\tilde{e}_i=\mathcal{U}\mathcal{D}(e_{i,\theta}):\mathbb{Z}\times X_{*}(T)\rightarrow X_{*}(T),
\]
\[
\tilde{\gamma}_i=\mathcal{U}\mathcal{D}(\gamma_i\circ \theta): X_{*}(T)\rightarrow \mathbb{Z},
\ \ 
\tilde{\varepsilon}_i=\mathcal{U}\mathcal{D}(\varepsilon_i\circ\theta): X_{*}(T)\rightarrow \mathbb{Z},
\]

\begin{thm}\label{UDthm}{\rm \cite{BK0, N0}}
Let $\chi=(X,\{e_i\}_{i\in [1,r]},\{\gamma_i\}_{i\in [1,r]},\{\varepsilon_i\}_{i\in [1,r]})$ be a geometric crystal,
$T$ an algebraic torus, $\theta:T\rightarrow X$ its positive structure. Then $(X_{*}(T), \{\tilde{e}_i\}_{i\in[1,r]},
\{\tilde{\gamma}_i\}_{i\in[1,r]},\{\tilde{\varepsilon}_i\}_{i\in[1,r]})$ has a free crystal structure.
\end{thm}

In the above notation, for $x\in X_*(T)$,
$x\mapsto \tilde{e}_i(1,x)$ and
$x\mapsto \tilde{e}_i(-1,x)$
give actions of Kashiwara operators
on $X_*(T)$, and
$\tilde{\gamma}_i$ define
the weight functions $x\mapsto \tilde{\gamma}_i(x)$ on $X_*(T)$.
The maps $\varphi_i$ are defined by $\varphi_i(x)=\tilde{\varepsilon}_i(x)+\tilde{\gamma}_i(x)$ $(i\in[1,r])$.

\subsection{Geometric crystal actions on $G^{u,e}$}\label{GCDB}

For a Weyl group element $u\in W$, we set $B_u^-:=B^-\cap N\ovl{u}N$. Note that
$G^{u,e}=B^-\cap B\ovl{u}B=HB_u^-$. Let $\gamma_i:G^{u,e}\rightarrow \mathbb{C}^{\times}$ be the rational function defined by 
\[
\gamma_i:G^{u,e} \hookrightarrow B^- \overset{\sim}{\rightarrow} H\times N^-
\overset{\rm proj}{\longrightarrow}H\overset{\al_i}{\rightarrow}\mathbb{C}^{\times}.
\]
For $\al\in \Delta^{\rm re}$, let $\mathfrak{g}_{\al}$ be the root space, and $N_{\al}:={\rm exp}(\mathfrak{g}_{\al})$.
For $i\in [1,r]$, we set $N_i^{\pm}:=N^{\pm}\cap \ovl{s}_iN^{\mp}\ovl{s}_i^{-1}$ and
$N^i_{\pm}:=N^{\pm}\cap \ovl{s}_iN^{\pm}\ovl{s}_i^{-1}$. Indeed, we have
$N^{\pm}_{i}=N_{\pm\al_i}$. We also set
\[
Y_{\pm\al_i}:=\lan x_{\pm i}(t) N_{\al} x_{\pm i}(-t) | t\in\mathbb{C}, \al\in\Delta^{\rm re}_{\pm}\setminus \{\pm\al_i\} \ran,
\]
where $x_{-i}(t):=y_i(t)$.
\begin{lem}{\rm \cite{Kum}}
For a simple root $\al_i$ $(i\in[1,r])$, we have

$(i)$ $Y_{\pm\al_i}=N^i_{\pm}$,

$(ii)$ $N^{\pm}=N_i^{\pm} \cdot Y_{\pm\al_i}$ (semi-direct product).
\end{lem}

By this lemma, we have the unique decomposition:
\[
N^-=N_i^- \cdot Y_{-\al_i}= N_{-\al_i} \cdot N^i_-,
\]
and we get the canonical projection $\xi_i:N^-\rightarrow N_{-\al_i}$.
Let $\chi_i$ be the function on $N^-$ defined as $\chi_i:=y_i^{-1}\circ \xi_i :N^- \rightarrow \mathbb{C}$,
where $y_i :\mathbb{C}\rightarrow N_{-\al_i}$ is defined as $c\mapsto y_i(c)$.
We extend this to the function on $B^-$ by $\chi_i(u\cdot t):=\chi_i(u)$ for $u\in N^-$ and $t\in H$.
We set
\[
\varphi_i:=(\chi_i|_{G^{u,e}})^{-1}: G^{u,e}\rightarrow \mathbb{C}^{\times},\quad 
\varepsilon_i:=\frac{\varphi_i}{\gamma_i}: G^{u,e}\rightarrow \mathbb{C}^{\times}.
\]
For a reduced expression $u=s_{i_1}\cdots s_{i_n}$, we suppose that $\{i_1,\cdots,i_n\}=\{1,2,\cdots,r\}$.
Then $\chi_i|_{G^{u,e}}$ is not identically zero \cite{N0}.
Thus, in this case, we can define the rational functions $\varphi_i$ and $\varepsilon_i$ $(i\in[1,r])$. 

For each $i\in[1,r]$,
let us define the rational $\mathbb{C}^{\times}$-action $e_i$ on $G^{u,e}$ as
($c\in\mathbb{C}^{\times}$, $x\in G^{u,e}$)
\[ e^c_i(x):=
x_i((c-1)\varphi_i(x)) x x_i((c^{-1}-1)\varepsilon_i(x)),
\]
if $\chi_i|_{G^{u,e}}\neq0$ and $e^c_i(x)=x$ if $\chi_i|_{G^{u,e}}=0$.
Let  $\al^{\vee}_{i}(T):=T^{\al^{\vee}_i}\in H$ for $T\in\mathbb{C}^{\times}$.

\begin{prop}
For $u\in W$ and its reduced expression $u=s_{i_1}\cdots s_{i_n}$,
the set
\[
\mathbb{B}_u^-:=\{t y_{i_1}(c_1)\al^{\vee}_{i_1}(c_1^{-1})\cdots 
y_{i_n}(c_n)\al^{\vee}_{i_n}(c_n^{-1})
|t\in H,\ c_1,\cdots,c_n\in\mathbb{C}^{\times}
\}
\]
is an open subset of $G^{u,e}$.
\end{prop}

\nd
{\it Proof.}

The set
\begin{equation}\label{anopen}
\ovl{\mathbb{B}}_u^-:=
\{t y_{i_1}(c_1)\cdots 
y_{i_n}(c_n)
|t\in H,\ c_1,\cdots,c_n\in\mathbb{C}^{\times}
\}
\end{equation}
is an open subset of $G^{u,e}$ \cite{HW}.
We prove our claim by showing $\ovl{\mathbb{B}}_u^-=\mathbb{B}_u^-$.
Note that
$y_{i}(S)\al^{\vee}_{j}(T^{-1})
=\al^{\vee}_{j}(T^{-1})y_{i}(ST^{-a_{j,i}})$
for $S$, $T\in\mathbb{C}^{\times}$ and $i$, $j\in\{1,2,\cdots,r\}$.
We have
\begin{eqnarray*}
& &
\mathbb{B}_u^-\ni t y_{i_1}(c_1)\al^{\vee}_{i_1}(c_1^{-1})\cdots 
y_{i_n}(c_n)\al^{\vee}_{i_n}(c_n^{-1})\\
&=& t \al^{\vee}_{i_1}(c_1^{-1})\cdots \al^{\vee}_{i_n}(c_n^{-1}) 
y_{i_1}(\frac{1}{c_1}\prod^{n}_{j=2}c_j^{-a_{j,1}})
\cdots y_{i_s}(\frac{1}{c_s}\prod^{n}_{j=s+1}c_j^{-a_{j,s}})
\cdots y_{i_n}(\frac{1}{c_n})\in \ovl{\mathbb{B}}_u^-,
\end{eqnarray*}
which means $\mathbb{B}_u^-\subset \ovl{\mathbb{B}}_u^-$.
Next, for $c_1,\cdots,c_n\in\mathbb{C}^{\times}$,
we put $\zeta_n:=c_n$ and
$\zeta_s:=\zeta_n^{-a_{i_n,i_{s}}}\zeta_{n-1}^{-a_{i_{n-1},i_{s}}}
\cdots \zeta_{s+1}^{-a_{i_{s+1},i_{s}}}c_s$ for $s=1,\cdots,n-1$.
It is clear that $\zeta_j\in\mathbb{C}^{\times}$ for $j=1,2,\cdots,n$.
We obtain
\begin{eqnarray*}
& &
\ovl{\mathbb{B}}_u^-\ni t y_{i_1}(c_1)\cdots 
y_{i_n}(c_n)\\
&=& t 
\al^{\vee}_{i_1}(\zeta_1^{-1})\cdots \al^{\vee}_{i_n}(\zeta_n^{-1}) 
\al^{\vee}_{i_1}(\zeta_1)\cdots \al^{\vee}_{i_n}(\zeta_n) 
y_{i_1}(c_1)\cdots 
y_{i_n}(c_n)\\
&=&t \al^{\vee}_{i_1}(\zeta_1^{-1})\cdots \al^{\vee}_{i_n}(\zeta_n^{-1})
y_{i_1}(\frac{c_1}{\zeta_1^2}\prod^{n}_{j=2}\zeta_{j}^{-a_{i_j,i_1}})
\al^{\vee}_{i_1}(\zeta_1)\\
& &
\cdots
y_{i_s}(\frac{c_s}{\zeta_s^2}\prod^{n}_{j=s+1}\zeta_{j}^{-a_{i_j,i_s}})
\al^{\vee}_{i_s}(\zeta_s) \cdots y_{i_n}(\zeta_n^{-1})
\al^{\vee}_{i_n}(\zeta_n)\\
&=&t \al^{\vee}_{i_1}(\zeta_1^{-1})\cdots \al^{\vee}_{i_n}(\zeta_n^{-1})
y_{i_1}(\zeta_1^{-1})\al^{\vee}_{i_1}(\zeta_1)
\cdots
y_{i_s}(\zeta_s^{-1})\al^{\vee}_{i_s}(\zeta_s)
\cdots
y_{i_n}(\zeta_n^{-1})\al^{\vee}_{i_n}(\zeta_n)\in \mathbb{B}_u^-,
\end{eqnarray*}
which means $\ovl{\mathbb{B}}_u^-\subset \mathbb{B}_u^-$. \qed

\begin{prop}{\rm \cite{N0,N}}\label{GCcell}
For $u\in W$ and its reduced expression $u=s_{i_1}\cdots s_{i_n}$, we suppose that $\{i_1,\cdots,i_n\}=\{1,2,\cdots,r\}$.
Then the quintuple $(G^{u,e},\{e_i\}_{i\in [1,r]},\{\gamma_i\}_{i\in [1,r]},\{\varepsilon_i\}_{i\in [1,r]})$ is a
geometric crystal.
The map $H\times (\mathbb{C}^{\times})^n\rightarrow G^{u,e}$,
$(t,c_1,\cdots,c_n)\mapsto t y_{i_1}(c_1)\al^{\vee}_{i_1}(c_1^{-1})\cdots y_{i_n}(c_n)\al^{\vee}_{i_n}(c_n^{-1})$
is a positive structure on this geometric crystal.
\end{prop}

\begin{prop}{\rm \cite{N0,N}}\label{geom-prop}
Let $u=s_{i_1}\cdots s_{i_n}$ be a reduced expression of $u\in W$ such that $\{i_1,\cdots,i_n\}=\{1,2,\cdots,r\}$.
The action of  $e^c_j$ on the open subset 
\[
\mathbb{B}_u^-=\{t y_{i_1}(t_1)\al^{\vee}_{i_1}(t_1^{-1}) y_{i_2}(t_2)\al^{\vee}_{i_2}(t_2^{-1})
\cdots y_{i_n}(t_n)\al^{\vee}_{i_n}(t_n^{-1})|t\in H,\ t_1,\cdots,t_n\in\mathbb{C}^{\times}\}\subset G^{u,e}
\]
is given by
\begin{eqnarray*}
& &e^c_j (t y_{i_1}(t_1)\al^{\vee}_{i_1}(t_1^{-1}) y_{i_2}(t_2)\al^{\vee}_{i_2}(t_2^{-1})\cdots y_{i_n}(t_n)\al^{\vee}_{i_n}(t_n^{-1}))\\
&=&t y_{i_1}(t_1')\al^{\vee}_{i_1}(t_1'^{-1}) y_{i_2}(t_2')\al^{\vee}_{i_2}(t_2'^{-1})\cdots y_{i_n}(t_n')\al^{\vee}_{i_n}(t_n'^{-1}),
\end{eqnarray*}
where 
\[
t_k'=t_k\frac{c \sum_{1\leq m<k, i_m=j}t_1^{a_{i_1,i}}\cdots t_{m-1}^{a_{i_{m-1},i}}t_m+\sum_{k\leq m\leq n, i_m=j}t_1^{a_{i_1,i}}\cdots t_{m-1}^{a_{i_{m-1},i}}t_m}{c \sum_{1\leq m\leq k, i_m=j}t_1^{a_{i_1,i}}\cdots t_{m-1}^{a_{i_{m-1},i}}t_m+\sum_{k< m\leq n, i_m=j}t_1^{a_{i_1,i}}\cdots t_{m-1}^{a_{i_{m-1},i}}t_m}.
\]
Furthermore,
\[
\varepsilon_j( (t y_{i_1}(t_1)\al^{\vee}_{i_1}(t_1^{-1})\cdots y_{i_n}(t_n)\al^{\vee}_{i_n}(t_n^{-1})))
=\left(\sum_{1\leq m\leq n,\ i_m=j}\frac{1}{t_m t^{a_{i_{m+1},j}}_{m+1}\cdots t^{a_{i_{n},j}}_{n}}\right)^{-1},
\]
\[
\gamma_j( (t y_{i_1}(t_1)\al^{\vee}_{i_1}(t_1^{-1})\cdots y_{i_n}(t_n)\al^{\vee}_{i_n}(t_n^{-1})))
=\frac{\al_j(t)}{t_1^{a_{i_1,j}}\cdots t_n^{a_{i_n,j}}}.
\]
\end{prop}


\subsection{Generalized minors and a bilinear form}\label{bilin}

We set $G_0:=N^-HN$, and let $x=[x]_-[x]_0[x]_+$ with $[x]_-\in N^-$, $[x]_0\in H$, $[x]_+\in N$ be the corresponding decomposition. 

\begin{defn}
For $i\in\{1,\cdots,\tilde{r}\}$ and $w,\ w'\in W$, the {\it generalized minor}
$\Delta_{w'\Lambda_i,w\Lambda_i}$ is a regular function on $G$
whose restriction to the open set $\ovl{w'}G_0\ovl{w}^{-1}$ is given by
$\Delta_{w'\Lambda_i,w\Lambda_i}(x)=([\ovl{w'}^{-1}x \ovl{w} ]_0)^{\Lambda_i}$.
Here, $\Lambda_i$ is the $i$-th  fundamental weight and for $a=T^h\in H$
($h\in \oplus_{j\in\{1,\cdots,\tilde{r}\}}\mathbb{Z}\al_j^{\vee}$, $T\in\mathbb{C}^{\times}$), we set $a^{\Lm_i}:=T^{\Lm_i(h)}$.

\end{defn}
Let $\omega:\ge\to\ge$ be the anti-involution 
\[
\omega(e_i)=f_i,\q
\omega(f_i)=e_i,\q \omega(h)=h,
\] and extend it to $G$ by setting
$\omega(x_i(c))=y_{i}(c)$, $\omega(y_{i}(c))=x_i(c)$ and $\omega(t)=t$
$(t\in H)$. One can calculate the generalized minors as follows. There exists a $\ge$ (or $G$)-invariant bilinear form on the
irreducible highest weight
$\ge$-module $V(\lm)$ such that 
\[
 \lan au,v\ran=\lan u,\omega(a)v\ran,
\q\q(u,v\in V(\lm),\,\, a\in \ge\ (\text{or }G)).
\]
For $g\in G$, 
we have the following simple fact:
\[
 \Del_{\Lm_i,\Lm_i}(g)=\lan gv_{\Lm_i},v_{\Lm_i}\ran,
\]
where $v_{\Lm_i}$ is a properly normalized highest weight vector in
$V(\Lm_i)$. Hence, for $w,\ w'\in W$, we have
\begin{equation}\label{minor-bilin}
 \Del_{w'\Lm_i,w\Lm_i}(g)=
\Del_{\Lm_i}(\ovl{w'}^{-1}g\ovl w)=
\lan g\ovl w\cdot v_{\Lm_i},\ovl{w'}\cdot v_{\Lm_i}\ran.
\end{equation}

Using generalized minors, the rational functions $\gamma_i$ and  $\varphi_i$ in \ref{GCDB}
are written as
\begin{equation}\label{minorexp}
\gamma_i=\prod^{\tilde{r}}_{j=1}\Delta_{\Lm_j,\Lm_j}^{a_{j,i}},\ \ 
\varphi_i=\frac{\Delta_{\Lm_i,\Lm_i}}{\Delta_{s_i\Lm_i,\Lm_i}},
\end{equation}
where we use (\ref{al-rel}) in the first relation.

\subsection{Glued crystal}\label{glued-cry-sec}

\begin{defn}\label{glued-cry}
Let $\{B_k\}_{k\in K}$ be a family of crystals,
$\{m_{k,k'}\}_{k,k'\in K}$ be a family of crystal isomorphisms 
$m_{k,k'}:B_k\overset{\sim}{\rightarrow} B_{k'}$ ($k,k'\in K$), where $K$ is an index set.
Then we call the set
\[
\mathcal{B}:=\amalg_{k\in K}B_k / \{{\rm identifications}\ m_{k,k'}\}
\]
 {\it glued crystal}.
\end{defn}

\begin{ex}

Let $\frak{g}$ be a symmetrizable Kac-Moody Lie algebra over $\mathbb{C}$ and 
$U_q(\frak{g})$ the associated quantum enveloping algebra.
Let ($B(\infty),\{\tilde{e}_i\},\{\tilde{f}_i\},\{\varepsilon_i\},\{\varphi_i\}, \{{\rm wt}_i\})$ be the crystal
associated with the crystal base of $U_q^-(\frak{g})$ and
$*:U_q(\frak{g})\rightarrow U_q(\frak{g})$ the antiautomorphism
such that $e_i^*=e_i$, $f_i^*=f_i$ and $(q^h)^*=q^{-h}$ in \cite{K}.
It is known that the map $*$ induces a bijection $*:B(\infty)\rightarrow B(\infty)$ satisfying $*\circ *=id$.
Let $B(\infty)^{*}$ be the crystal as follows : $B(\infty)^{*}$ is equal to $B(\infty)$ as sets,
and maps are defined as $\tilde{e}_i^*:=*\circ \tilde{e}_i\circ *$, 
$\tilde{f}_i^*:=*\circ \tilde{f}_i\circ *$,
$\varepsilon_i^*:=\varepsilon_i\circ *$,
$\varphi_i^*:=\varphi_i\circ *$ and ${\rm wt}_i^*:={\rm wt}_i\circ *$.
Clearly, $*:B(\infty)\rightarrow B(\infty)^*$ is a crystal isomorphism. Thus we get
a glued crystal
\[
B(\infty)\sqcup B(\infty)^*/({\rm identification\ by}\ *).
\]


\end{ex}

\section{Cluster ensembles}\label{ensect}

Following \cite{HW}, let us recall the notion of the cluster $\mathcal{A}$-variety and $\mathcal{X}$-variety.

\subsection{Definitions of cluster $\mathcal{A}$-variety and $\mathcal{X}$-variety }

\begin{defn}
A {\it seed} $\Sigma=(I,I_0,B,d)$ is a quintuple of the following data:
\begin{enumerate}
\item $I$ is a finite index set and $I_0$ is a subset of $I$. The elements of $I_0$ are called frozen.
\item $B=(b_{i,j})$ is an $I\times I$-matrix called {\it exchange matrix} which satisfies $b_{i,j}\in\mathbb{Z}$ unless both $i$ and $j$ are frozen.
\item $d=(d_i)_{i\in I}$ is a set of positive integers such that $b_{i,j}d_j=-b_{j,i}d_i$.
\end{enumerate}
\end{defn}

For $k\in I\setminus I_0$, we say a seed $\Sigma'=(I',I_0,B',d')$ is obtained from $\Sigma$ by {\it mutation} at $k$ if there exists a bijective
$M_k:I\rightarrow I'$ which satisfies $M_k (i)=i$ for $i\in I_0$, $d'_{M_k(j)}=d_j$ for $j\in I$ and 
\[b_{M_k(i),M_k(j)}':=
\begin{cases}
	-b_{ij} & {\rm if}\ i=k\ {\rm or}\ j=k, \\
	b_{ij}+\frac{|b_{ik}|b_{kj}+b_{ik}|b_{kj}|}{2} & {\rm otherwise}.
\end{cases}
\]
Then we write $\Sigma'=\mu_k(\Sigma)$.

To a seed $\Sigma=(I,I_0,B,d)$, we associate a set of variables $\{A_i\}_{i\in I}$ and a split algebraic torus $\mathcal{A}_{\Sigma}={\rm Spec}\ \mathbb{C}[A_i^{\pm1}|i\in I]$ called {\it cluster $\mathcal{A}$-torus }. We call $\{A_i\}_{i\in I}$ {\it cluster $\mathcal{A}$-coordinate} on $\mathcal{A}_{\Sigma}$. If $\Sigma'=(I',I_0,B',d')=\mu_k(\Sigma)$ then there exists a birational map (called {\it mutation}) $\mu_k:\mathcal{A}_{\Sigma}\rightarrow \mathcal{A}_{\Sigma'}$ defined as
\[
\mu_k^*(A'_{M_k(i)})=
\begin{cases}
A_i & {\rm if}\ i\neq k,\\
\frac{\prod_{b_{k,j}>0}A_j^{b_{k,j}}+\prod_{b_{k,j}<0}A_j^{-b_{k,j}}}{A_k} & {\rm if}\ i=k,
\end{cases}
\]
where $\{A'_i\}_{i\in I'}$ is the cluster $\mathcal{A}$-coordinate associated to $\Sigma'$. 

\begin{defn}
Cluster $\mathcal{A}$-variety $\mathcal{A}_{|\Sigma|}$ is the scheme obtained by gluing together
all tori $\mathcal{A}_{\Sigma'}$ of seeds $\Sigma'$ which are obtained from $\Sigma$ by an iteration of mutations.
\end{defn}

To a seed $\Sigma=(I,I_0,B,d)$, we also associate an algebraic torus
$\mathcal{X}_{\Sigma}:={\rm Spec}\ \mathbb{C}[X_i^{\pm1}|i\in I]$
called {\it cluster $\mathcal{X}$-torus} with variables $\{X_i\}_{i\in I}$.
We call $\{X_i\}_{i\in I}$ {\it cluster $\mathcal{X}$-coordinate} on $\mathcal{X}_{\Sigma}$.
If $\Sigma'=\mu_k(\Sigma)$ then there exists a birational map (called {\it mutation})
 $\mu_k:\mathcal{X}_{\Sigma}\rightarrow \mathcal{X}_{\Sigma'}$ defined as
\[
\mu_k^*(X'_{M_k(i)})=
\begin{cases}
X_iX_k^{[b_{i,k}]_+}(1+X_k)^{-b_{i,k}} & {\rm if}\ i\neq k,\\
X_k^{-1} & {\rm if}\ i=k,
\end{cases}
\]
where $\{X'_i\}_{i\in I'}$ is the cluster $\mathcal{X}$-coordinate associated to $\Sigma'$ and $[b_{i,k}]_+:={\rm max} (b_{i,k},0)$.

\begin{defn}
Cluster $\mathcal{X}$-variety $\mathcal{X}_{|\Sigma|}$ is the scheme obtained by gluing together
all tori $\mathcal{X}_{\Sigma'}$ of seeds $\Sigma'$ which are obtained from $\Sigma$
by an iteration of mutations.
\end{defn}

In what follows, we identify $I$ with $I'$ by $M_k$, and write
$b_{M_k(i),M_k(j)}'=b_{i,j}'$, $A'_{M_k(i)}=A'_i$ and $X'_{M_k(i)}=X'_i$.

\begin{prop}\label{map-p}{\rm \cite{HW}}
Let $M=(M_{i,j})$ be an $I\times I$-matrix such that $M_{i,j}=0$ unless both $i$ and $j$ are frozen.
For a seed $\Sigma=(I,I_0,B,d)$ such that $\tilde{B}=B+M$ is an integer matrix,
we define a map $p_M:\mathcal{A}_{\Sigma}\rightarrow \mathcal{X}_{\Sigma}$ as
\[
p_M^* (X_i)=\prod_{j\in I}A^{\tilde{B}_{i,j}}_j.
\]
Then $p_M$ extends to a regular map $p_M:\mathcal{A}_{|\Sigma|}\rightarrow \mathcal{X}_{|\Sigma|}$.
\end{prop}

If $\Sigma'=\mu_k(\Sigma)$ and $B'$ is the exchange matrix of $\Sigma'$
then $\tilde{B'}=B'+M'$ is an integer matrix.
Thus, we can define
$p'_M:\mathcal{A}_{\Sigma'}\rightarrow \mathcal{X}_{\Sigma'}$ by $p_M'^* (X'_i)=\prod_{j\in I}A'^{\tilde{B}'_{i,j}}_j$.
This proposition means the following diagram is commutative:
\begin{equation}\label{commp}
  \xymatrix{
    \mathcal{A}_{\Sigma} \ar[r]^{\mu_k} \ar[d]_{p_M} & \mathcal{A}_{\Sigma'} \ar[d]^{p_M'}  \\
    \mathcal{X}_{\Sigma} \ar[r]_{\mu_k} & \mathcal{X}_{\Sigma'}
  }
\end{equation}

\nd
The map $p_M:\mathcal{A}_{|\Sigma|}\rightarrow \mathcal{X}_{|\Sigma|}$ is called an {\it ensemble map} in the context of \cite{HW}.

\subsection{The set of $\mathbb{Z}^T$-valued points}\label{zvalue}

Let $\mathbb{Z}^T$ be the tropical semi-field of integers, that is,
it is equal to $\mathbb{Z}$ as sets, and product and sum are defined as
$+$ and {\rm max} respectively.
For a split torus $H$, we set
\[
H(\mathbb{Z}^T):=X_*(H),
\]
where $X_*(H)$ is the group of co-characters of $H$.
We can verify that $H(\mathbb{Z}^T)\cong (\mathbb{Z}^T)^{{\rm dim} H}$.
Note that a positive rational map $f:H\rightarrow H'$ induces a map
$\mathcal{U}\mathcal{D}(f):H(\mathbb{Z}^T)\rightarrow H'(\mathbb{Z}^T)$,
where $\mathcal{U}\mathcal{D}$ is the functor in \ref{GCsub}.

\begin{defn}\cite{FG}
For a seed $\Sigma$ and the cluster $\mathcal{A}$-variety $\mathcal{A}_{|\Sigma|}$ and
$\mathcal{X}$-variety $\mathcal{X}_{|\Sigma|}$, we set
\[
\mathcal{A}_{|\Sigma|}(\mathbb{Z}^T):=
\coprod \mathcal{A}_{\Sigma'}(\mathbb{Z}^T)/\{{\rm identifications}\ \mathcal{U}\mathcal{D}(\mu),\ \mu:{\rm mutation}\} ,
\]
\[
\mathcal{X}_{|\Sigma|}(\mathbb{Z}^T):=
\coprod \mathcal{X}_{\Sigma'}(\mathbb{Z}^T)/\{{\rm identifications}\ \mathcal{U}\mathcal{D}(\mu),\ \mu:{\rm mutation}\} ,
\]
where $\Sigma'$ runs over the set of all seeds which are obtained from $\Sigma$
by an iteration of mutations. We call them the sets of  $\mathbb{Z}^T$-{\it valued points} of
cluster $\mathcal{A}$($\mathcal{X}$)-variety.
\end{defn}

The sets of  $\mathbb{Z}^T$-valued points of
cluster varieties appear in the context of 
Fock-Goncharov conjecture in \cite{FG,GHKK}.

\subsection{Seeds associated with reduced words}\label{seedass}

We will use the notation in \ref{notation}. For $u\in W$ and its reduced word $\textbf{i}=(i_1,\cdots,i_n)$, one associate a seed $\Sigma_{\textbf{i}}$ as follows. We set $i_{-j}=-j$ for $j=1,2,\cdots,\tilde{r}$.

\begin{defn}\label{exmatdef}\cite{HW}
We define the index set as $I:=\{-\tilde{r},\cdots,-2,-1\}\cup\{1,2,\cdots,n\}$. For $k\in I$, we set $k^+:={\rm min}\{l\in I | l>k, |i_l|=|i_k|\}\cup\{n+1\}$ and $I_0:=\{k\in I| k<0,\ {\rm or}\ k^+>n\}$. The exchange matrix $B_{\textbf{i}}=(b_{j,k})$ is defined by
\begin{eqnarray*}
b_{j,k}&=&\frac{a_{|i_k|,|i_j|}}{2} (-[j=k^+]+[j^+=k]-[k<j<k^+][j>0]+[k<j^+<k^+][j^+\leq n]\\
& &+[j<k<j^+][k>0]-[j<k^+<j^+][k^+\leq n]),
\end{eqnarray*}
where for a proposition $P$, 
\[
[P]=
\begin{cases}
1 & {\rm if}\ P:{\rm true},\\
0 & {\rm if}\ P:{\rm false}.
\end{cases}
\]
Let $d_k=d_{|i_k|}$ for $k\in I$, where 
the right-hand side means the symmetrizer of the Cartan matrix $(a_{i,j})_{i,j\in I}$.
Then we define a seed $\Sigma_{\textbf{i}}:=(I,I_0,B_{\textbf{i}},d)$.
\end{defn}

The lattice $\bigoplus_{1\leq i\leq \tilde{r}}\mathbb{Z}\al_i$
is a sublattice of $P$, and we see that 
its kernel $\{t\in H| t^{\al_i}=1,\ 1\leq i\leq \tilde{r}\}$
is a discrete subgroup of the center of $G$.
Let $G_{{\rm Ad}}$ denote the quotient of $G$ by the discrete subgroup
$\{t\in H| t^{\al_i}=1,\ 1\leq i\leq \tilde{r}\}$ and $H_{\rm Ad}$ denote the image of $H$ in $G_{\rm Ad}$.
The character lattice ${\rm Hom}(H_{\rm Ad},\mathbb{C}^{\times})$ is
canonically isomorphic to $\bigoplus_{1\leq i\leq \tilde{r}}\mathbb{Z}\al_i$.
Thus, the co-character lattice ${\rm Hom}(\mathbb{C}^{\times},H_{\rm Ad})$ has
a dual basis $\Lm_1^{\vee},\cdots,\Lm_{\tilde{r}}^{\vee}$ of
 {\it fundamental coweights} 
such that
$\al_i(\Lm_j^{\vee})=\delta_{i,j}$
for $i,j \in[1,\tilde{r}]$. 
Let $T^{\Lm_i^{\vee}}$ be 
an element
of $H_{\rm Ad}$ such that $\al_{j}(T^{\Lm_i^{\vee}})=T^{\delta_{i,j}}$ for
$i,j\in [1,\tilde{r}]$ and $T\in\mathbb{C}^{\times}$.
Now we define numbers $a_{ij}$ $(i,j \in[1,\tilde{r}])$ as
$a_{i,j}:=\al_j(\al_i^{\vee})$.

\begin{prop}{\rm \cite{HW}}\label{ExDyn}
The $\tilde{r}\times \tilde{r}$ integer matrix $(a_{i,j})_{1\leq i,j\leq \tilde{r}}$ is nondegenerate and symmetriz-
able. We also get $\al_i^{\vee}=\sum^{\tilde{r}}_{j=1} a_{ij}\Lm_j^{\vee}$.
\end{prop}

\begin{defn}\cite{HW}\label{biratio-def1}
Let $\mathcal{X}_{\Sigma_{\textbf{i}}}$ be the cluster $\mathcal{X}$-torus which associates to 
the seed $\Sigma_{\textbf{i}}$.
An open immersion $x_{\Sigma_{\bf{{\rm i}}}}:\mathcal{X}_{\Sigma_{\textbf{i}}}\rightarrow G^{u,e}_{\rm Ad}$ is defined
for an element $u\in W$ and its reduced word $\textbf{i}=(i_1,i_2,\cdots,i_n)$ :
\[
x_{\Sigma_{\bf{{\rm i}}}}: (X_{-\tilde{r}},\cdots,X_{-1},X_1,\cdots,X_n)
\mapsto X_{-\tilde{r}}^{\Lm_{\tilde{r}}^{\vee}}\cdots X_{-1}^{\Lm_{1}^{\vee}}
y_{i_1}(1) X_1^{\Lm_{i_1}^{\vee}} y_{i_2}(1) X_2^{\Lm_{i_2}^{\vee}}\cdots y_{i_n}(1) X_n^{\Lm_{i_n}^{\vee}}.
\]
\end{defn}

Let $u_{\leq k}:=s_{i_1}\cdots s_{i_k}$ for $k\in[1,n]$ and $u_{\leq k}=e$ for $k\in\{-\tilde{r},\cdots,-2,-1\}$.

\begin{lem}{\rm \cite{HW}}\label{biratio-def2}
Let $\mathcal{A}_{\Sigma_{{\rm \bf{i}}}}$ be the cluster $\mathcal{A}$-torus which associates to 
the seed $\Sigma_{\textbf{i}}$.
There exists an open immersion $a_{\Sigma_{{\rm \bf{i}}}}:\mathcal{A}_{\Sigma_{{\rm \bf{i}}}}\rightarrow G^{u,e}$
such that the pull-back $a_{\Sigma_{{\rm \bf{i}}}}^*$ identifies each coordinate function $A_k$ with a
generalized minor $\Delta_{u_{\leq k}\Lm_{|i_k|},\Lm_{|i_k|}}$ for $k\in I$.
\end{lem}

\begin{defn}\cite{HW} For $x\in G_0:=N^-HN$, we write $x=[x]_- [x]_0[x]_+$ with $[x]_-\in N^-$, $[x]_0\in H$, $[x]_+\in N_+$.
For $u\in W$, the twist map $\zeta^{u,e}:G^{u,e}\rightarrow G^{u^{-1},e}$ is defined by
\[ x\mapsto \theta([\overline{u}^{-1}x]_-^{-1} \overline{u}^{-1}x), \]
where $\theta$ is the automorphism of $G$ such that $\theta(a)=a^{-1}$ $(a\in H)$, $\theta(x_{i}(T))=y_{i}(T)$ 
and $\theta(y_{i}(T))=x_{i}(T)$ $(T\in\mathbb{C})$.
We also define $\iota$ as the antiautomorphism of $G$ defined
by $a\mapsto a^{-1}$ for $a\in H$ and $x_{i}(T)\mapsto x_{i}(T)$, $y_{i}(T)\mapsto y_{i}(T)$
for $T\in\mathbb{C}$.
\end{defn}

\begin{prop}{\rm \cite{HW}}\label{twistbire}
For $x\in G^{u,e}$, we get $\ovl{u}^{-1}x\in G_0$.
The map $\iota\circ \zeta^{u,e}:G^{u,e}\rightarrow G^{u,e}$ is a biregular isomorphism.
\end{prop}

\begin{thm}\label{pgdef}{\rm \cite{HW}}
Let $M=(M_{j,k})$ be the following $I\times I$-matrix 
\[
M_{j,k}=\frac{a_{|i_k|,|i_j|}}{2} ([j^+,k^+>n]+[j,k<0]).
\]
\begin{enumerate}
\item[(1)] There is a regular map $a_{|\Sigma_{\bf{{\rm i}}}|}:\mathcal{A}_{|\Sigma_{\bf{{\rm i}}}|}\rightarrow G^{u,e}$
which extends $a_{\Sigma_{\bf{{\rm i}}}}:\mathcal{A}_{\Sigma_{\bf{{\rm i}}}}\rightarrow G^{u,e}$
in Lemma \ref{biratio-def2}. 
It induces an algebra isomorphism
$\mathbb{C}[G^{u,e}]\rightarrow \mathbb{C}[\mathcal{A}_{|\Sigma_{\bf{{\rm i}}}|}]$.
\item[(2)] There is a regular map $x_{|\Sigma_{\bf{{\rm i}}}|}:\mathcal{X}_{|\Sigma_{\bf{{\rm i}}}|}\rightarrow G^{u,e}_{{\rm Ad}}$
which extends $x_{\Sigma_{\bf{{\rm i}}}}:\mathcal{X}_{\Sigma_{\bf{{\rm i}}}}\rightarrow G^{u,e}_{{\rm Ad}}$
in Definition \ref{biratio-def1}.
\item[(3)] The all entries of $B_{{\rm \bf{i}}}+M$ are integer and
the following diagram is commutative:
\[
  \xymatrix{
    \mathcal{A}_{|\Sigma_{\bf{{\rm i}}}|} \ar[r]^{a_{|\Sigma_{\bf{{\rm i}}}|}} \ar[d]_{p_M} & G^{u,e} \ar[d]^{p_G}  \\
    \mathcal{X}_{|\Sigma_{\bf{{\rm i}}}|} \ar[r]_{x_{|\Sigma_{\bf{{\rm i}}}|}} & G^{u,e}_{{\rm Ad}}
  }
\]
where $p_G:G^{u,e}\rightarrow G^{u,e}_{{\rm Ad}}$ is the composition of
the automorphism $\iota\circ \zeta^{u,e}$ and the quotient map $G$ to $G_{{\rm Ad}}$.
\end{enumerate}
\end{thm}

\section{Compatibility between geometric crystal structures on the cluster tori}

In \ref{seedass}, we defined tori $\mathcal{A}_{\Sigma}$ (resp. $\mathcal{X}_{\Sigma}$)
birationally isomorphic to $G^{u,e}$
(resp. $G^{u,e}_{\rm Ad}$).
In this section, we will define two geometric crystal structures on each torus
$\mathcal{A}_{\Sigma}$ (resp. $\mathcal{X}_{\Sigma}$) (Definition \ref{actions}). 
Furthermore, we will discuss a compatibility between these geometric
crystal structures (Proposition \ref{compprop}).
In Sect. \ref{action-X}, \ref{action-A}, we will calculate
the explicit formulae of them. In the rest of article, 
for a reduced expression $u=s_{i_1}\cdots s_{i_n}$, 
we suppose that $\{i_1,\cdots,i_n\}=\{1,2,\cdots,r\}$.
Let $Z= \{t\in H| t^{\al_i}=1,\ 1\leq i\leq \tilde{r}\}$. 

First, we define a geometric crystal structure on $G^{u,e}_{\rm Ad}$.
Recall that $G^{u,e}$ has a geometric crystal structure
$(G^{u,e},\{e_j\}_{j\in[1,r]},\ \{\varepsilon_{j}\}_{j\in[1,r]},\ \{\gamma_j\}_{j\in[1,r]})$ (Proposition \ref{GCcell}).
The definitions imply that
$\varepsilon_{j}$ and $\gamma_j$ are $Z$-invariant, and
$e_j$ satisfies $e_j(c,tx)=t e_j(c,x)$ for $t\in Z$, $(c,x)\in {\rm dom}(e_j)$. 
Thus, the geometric crystal structure on $G^{u,e}$ induces a  
geometric crystal structure on $G^{u,e}_{\rm Ad}$.

\begin{defn}\label{GCdefX}.
We write this geometric crystal structure as
$(G^{u,e}_{\rm Ad},\{e_j\}_{j\in[1,r]},\ \{\varepsilon_{j}\}_{j\in[1,r]},\ \{\gamma_j\}_{j\in[1,r]})$,
that is, 
we use the same notation for $e_j$, $\varepsilon_{j}$ and $\gamma_j$ on $G^{u,e}_{\rm Ad}$ as
those on $G^{u,e}$.
\end{defn}

In the formulae of Proposition \ref{geom-prop} for $e_j^c$, $\varepsilon_j$ and $\gamma_j$,
an element $t\in H$ is replaced with $t\in H_{\rm Ad}$ when we consider the geometric crystal
structure on $G^{u,e}_{\rm Ad}$.





\begin{prop}
A biregular map $G^{u,e}_{\rm Ad}\rightarrow G^{u,e}_{\rm Ad}$ defined as
$x Z\mapsto (\iota\circ \zeta^{u,e}(x))Z$ $(x\in G^{u,e})$ is well-defined.
Let us denote it by $\overline{\iota\circ \zeta^{u,e}}$.
\end{prop}

\nd
{\it Proof.}

We take $x\in G^{u,e}$ and $z\in Z$. By Proposition \ref{twistbire}, we have $\ovl{u}^{-1}x\in G_0$.
Considering the decomposition $\ovl{u}^{-1}x=[\ovl{u}^{-1}x]_-[\ovl{u}^{-1}x]_0[\ovl{u}^{-1}x]_+$,
we obtain $\ovl{u}^{-1}xz=[\ovl{u}^{-1}x]_-[\ovl{u}^{-1}x]_0z[\ovl{u}^{-1}x]_+$ and
$[\ovl{u}^{-1}x]_0z\in H$, which yields
\begin{equation}\label{pr-0-1}
[\ovl{u}^{-1}xz]_-=[\ovl{u}^{-1}x]_-.
\end{equation}
It follows from (\ref{pr-0-1}) that
\begin{eqnarray*}
(\iota\circ\zeta^{u,e})(xz)&=&\iota \circ\theta([\ovl{u}^{-1}xz]^{-1}_{-}\ovl{u}^{-1}xz)\\
&=&\iota \circ\theta([\ovl{u}^{-1}x]^{-1}_{-}\ovl{u}^{-1}xz)\\
&=&\iota (\theta([\ovl{u}^{-1}x]^{-1}_{-}\ovl{u}^{-1}x)z^{-1})\\
&=&\iota(\theta([\ovl{u}^{-1}x]^{-1}_{-}\ovl{u}^{-1}x)) z\\
&=& (\iota\circ\zeta^{u,e})(x) z.
\end{eqnarray*}


\nd
Thus, the biregular map $G^{u,e}_{\rm Ad}\rightarrow G^{u,e}_{\rm Ad}$ is induced. \qed




\begin{defn}\label{actions}
We use the same notation as in Theorem \ref{pgdef}. 
Let $\Sigma$ be a seed obtained from $\Sigma_{{\rm \bf{i}}}$ by an iteration
of mutations $\ovl{\mu}$,
and $\ovl{\mu}^a:\mathcal{A}_{\Sigma_{\bf{{\rm i}}}}\rightarrow \mathcal{A}_{\Sigma}$,
 $\ovl{\mu}^x:\mathcal{X}_{\Sigma_{\bf{{\rm i}}}}\rightarrow \mathcal{X}_{\Sigma}$
be the corresponding birational maps. 
Let $a_{\Sigma}$, $x_{\Sigma}$ denote the birational maps
$a_{\Sigma_{\bf{{\rm i}}}}\circ(\ovl{\mu}^a)^{-1}:\mathcal{A}_{\Sigma}\rightarrow G^{u,e}$,
$x_{\Sigma_{\bf{{\rm i}}}}\circ(\ovl{\mu}^x)^{-1}:\mathcal{X}_{\Sigma}\rightarrow G^{u,e}_{\rm Ad}$,
respectively. We define two geometric crystal structures on the torus $\mathcal{A}_{\Sigma}$ (resp. $\mathcal{X}_{\Sigma}$)
as follows:

(1) The first one is 
\[
(\mathcal{A}_{\Sigma},\ a_{\Sigma}^{-1}\circ e_j^c\circ a_{\Sigma},
\ \varepsilon_j\circ a_{\Sigma},
\ \gamma_j\circ a_{\Sigma}),\quad
{\rm(resp.}\  
(\mathcal{X}_{\Sigma},\ x_{\Sigma}^{-1}\circ e_j^c\circ x_{\Sigma},
\ \varepsilon_j\circ x_{\Sigma},
\ \gamma_j\circ x_{\Sigma})). \] 
(2) The second one is
\[
(\mathcal{A}_{\Sigma},\ a_{\Sigma}^{-1}\circ(\iota\circ \zeta^{u,e})^{-1}\circ e_j^c
\circ (\iota\circ \zeta^{u,e})\circ a_{\Sigma},
\ \varepsilon_j\circ (\iota\circ \zeta^{u,e})\circ a_{\Sigma},
\ \gamma_j\circ (\iota\circ \zeta^{u,e})\circ a_{\Sigma}),
\]
\[
{\rm (resp.}\ 
(\mathcal{X}_{\Sigma},\ x_{\Sigma}^{-1}\circ (\ovl{\iota\circ \zeta^{u,e}})\circ 
e_j^c \circ (\ovl{\iota\circ \zeta^{u,e}})^{-1}\circ x_{\Sigma},\ \varepsilon_j\circ 
(\ovl{\iota\circ \zeta^{u,e}})^{-1}\circ x_{\Sigma},\ \gamma_j\circ (\ovl{\iota\circ \zeta^{u,e}})^{-1}\circ x_{\Sigma})).
\] 
\end{defn}


\begin{prop}\label{compprop}
Let $p=p_{\Sigma}$ be the restriction of the map $p_M$ in Theorem \ref{pgdef} to the torus $\mathcal{A}_{\Sigma}$.
The following commutative diagrams hold:
\begin{equation}\label{comp1}
\begin{xy}
(0,0)*{\mathcal{X}_{\Sigma}}="c",
(0,20) *{ \mathcal{A}_{\Sigma}}="a",
(60,0)*{\mathcal{X}_{\Sigma}}="d",
(60,20) *{ \mathcal{A}_{\Sigma}}="b",
\ar@{->} "a";"b"^{a_{\Sigma}^{-1}\circ (\iota\circ \zeta^{u,e})^{-1} \circ e_j^c \circ \iota\circ \zeta^{u,e} \circ a_{\Sigma}}
\ar@{->} "a";"c"_{p}
\ar@{->} "b";"d"^{p}
\ar@{->} "c";"d"_{x_{\Sigma}^{-1}\circ e_j^c \circ x_{\Sigma}}
\end{xy}\qquad
\begin{xy}
(0,0)*{\mathcal{X}_{\Sigma}}="c",
(0,20) *{ \mathcal{A}_{\Sigma}}="a",
(60,0)*{\mathcal{X}_{\Sigma}}="d",
(60,20) *{ \mathcal{A}_{\Sigma}}="b",
\ar@{->} "a";"b"^{a_{\Sigma}^{-1}\circ e_j^c \circ a_{\Sigma}}
\ar@{->} "a";"c"_{p}
\ar@{->} "b";"d"^{p}
\ar@{->} "c";"d"_{x_{\Sigma}^{-1}\circ (\ovl{\iota\circ \zeta^{u,e}}) \circ e_j^c \circ (\ovl{\iota\circ \zeta^{u,e}})^{-1} 
\circ x_{\Sigma}}
\end{xy}
\end{equation}

We also get the following commutative diagrams:
\begin{equation}\label{comp2}
\begin{xy}
(0,0)*{\mathcal{X}_{\Sigma}}="c",
(0,20) *{ \mathcal{A}_{\Sigma}}="a",
(30,0)*{\mathbb{C}^{\times}}="d",
\ar@{->} "a";"c"_{p}
\ar@{->} "a";"d"^{ \varepsilon_j\circ(\iota\circ \zeta^{u,e})\circ a_{\Sigma}}
\ar@{->} "c";"d"_{\varepsilon_j\circ x_{\Sigma}}
\end{xy}\qquad
\begin{xy}
(0,0)*{\mathcal{X}_{\Sigma}}="c",
(0,20) *{ \mathcal{A}_{\Sigma}}="a",
(30,0)*{\mathbb{C}^{\times}}="d",
\ar@{->} "a";"c"_{p}
\ar@{->} "a";"d"^{ \gamma_j\circ(\iota\circ \zeta^{u,e})\circ a_{\Sigma}}
\ar@{->} "c";"d"_{\gamma_j\circ x_{\Sigma}}
\end{xy}
\end{equation}
\[
\begin{xy}
(0,0)*{\mathcal{X}_{\Sigma}}="c",
(0,20) *{ \mathcal{A}_{\Sigma}}="a",
(30,0)*{\mathbb{C}^{\times}}="d",
\ar@{->} "a";"c"_{p}
\ar@{->} "a";"d"^{ \varepsilon_j\circ a_{\Sigma}}
\ar@{->} "c";"d"_{\varepsilon_j\circ(\ovl{\iota\circ\zeta^{u,e}})^{-1}\circ x_{\Sigma}}
\end{xy}\qquad
\begin{xy}
(0,0)*{\mathcal{X}_{\Sigma}}="c",
(0,20) *{ \mathcal{A}_{\Sigma}}="a",
(30,0)*{\mathbb{C}^{\times}}="d",
\ar@{->} "a";"c"_{p}
\ar@{->} "a";"d"^{ \gamma_j\circ a_{\Sigma}}
\ar@{->} "c";"d"_{\gamma_j\circ(\ovl{\iota\circ\zeta^{u,e}})^{-1}\circ x_{\Sigma}}
\end{xy}
\]
\end{prop}

\nd
{\it Proof.}

First, let us prove them for $\Sigma=\Sigma_{\bf{{\rm i}}}$. We set $p:=p_{\Sigma_{\bf{{\rm i}}}}$.
It follows from Theorem \ref{pgdef} (3), definitions of $p_G$ and action of $e_j^c$ that
\begin{eqnarray*}
x_{\Sigma_{\bf{{\rm i}}}}^{-1}\circ e_j^c \circ x_{\Sigma_{\bf{{\rm i}}}} \circ p
&=&x_{\Sigma_{\bf{{\rm i}}}}^{-1}\circ e_j^c \circ p_G\circ a_{\Sigma_{\bf{{\rm i}}}} \\
&=&x_{\Sigma_{\bf{{\rm i}}}}^{-1}\circ e_j^c \circ q \circ (\iota\circ \zeta^{u,e})\circ a_{\Sigma_{\bf{{\rm i}}}} \\
&=& x_{\Sigma_{\bf{{\rm i}}}}^{-1}\circ q\circ e_j^c \circ (\iota\circ \zeta^{u,e})\circ a_{\Sigma_{\bf{{\rm i}}}} \\
&=& x_{\Sigma_{\bf{{\rm i}}}}^{-1}\circ q\circ (\iota\circ \zeta^{u,e})\circ
 (\iota\circ \zeta^{u,e})^{-1} \circ e_j^c \circ (\iota\circ \zeta^{u,e})\circ a_{\Sigma_{\bf{{\rm i}}}} \\
&=& x_{\Sigma_{\bf{{\rm i}}}}^{-1}\circ p_G\circ (\iota\circ \zeta^{u,e})^{-1} \circ 
e_j^c \circ (\iota\circ \zeta^{u,e})\circ a_{\Sigma_{\bf{{\rm i}}}} \\
&=& p\circ a_{\Sigma_{\bf{{\rm i}}}}^{-1}\circ(\iota\circ \zeta^{u,e})^{-1} \circ 
e_j^c \circ (\iota\circ \zeta^{u,e})\circ a_{\Sigma_{\bf{{\rm i}}}},
\end{eqnarray*}
and
\begin{eqnarray*}
p\circ a_{\Sigma_{\bf{{\rm i}}}}^{-1}\circ e_j^c \circ a_{\Sigma_{\bf{{\rm i}}}}
&=&x_{\Sigma_{\bf{{\rm i}}}}^{-1}\circ p_G\circ e_j^c \circ a_{\Sigma_{\bf{{\rm i}}}} \\
&=&x_{\Sigma_{\bf{{\rm i}}}}^{-1}\circ (\ovl{\iota\circ \zeta^{u,e}}) \circ q\circ e_j^c 
\circ a_{\Sigma_{\bf{{\rm i}}}} \\
&=&x_{\Sigma_{\bf{{\rm i}}}}^{-1}\circ (\ovl{\iota\circ \zeta^{u,e}}) \circ e_j^c \circ 
q\circ a_{\Sigma_{\bf{{\rm i}}}} \\
&=&x_{\Sigma_{\bf{{\rm i}}}}^{-1}\circ (\ovl{\iota\circ \zeta^{u,e}}) \circ e_j^c \circ 
(\ovl{\iota\circ \zeta^{u,e}})^{-1}\circ (\ovl{\iota\circ \zeta^{u,e}})\circ
q \circ a_{\Sigma_{\bf{{\rm i}}}} \\
&=&x_{\Sigma_{\bf{{\rm i}}}}^{-1}\circ (\ovl{\iota\circ \zeta^{u,e}}) \circ e_j^c \circ 
(\ovl{\iota\circ \zeta^{u,e}})^{-1}\circ p_G \circ a_{\Sigma_{\bf{{\rm i}}}} \\
&=&x_{\Sigma_{\bf{{\rm i}}}}^{-1}\circ (\ovl{\iota\circ \zeta^{u,e}}) \circ e_j^c \circ 
(\ovl{\iota\circ \zeta^{u,e}})^{-1}\circ x_{\Sigma_{\bf{{\rm i}}}}\circ p.
\end{eqnarray*}
Thus, we get (\ref{comp1}) for $\Sigma=\Sigma_{\bf{{\rm i}}}$.

The definitions of $\gamma_j$ on $G^{u,e}$ and $G^{u,e}_{\rm Ad}$ mean
$\gamma_j\circ q (g)=\gamma_j(g)$ for $g\in G^{u,e}$.
 Hence, we have
\begin{eqnarray*}
\gamma_j\circ x_{\Sigma_{\rm\bf{i}}}\circ p &=& \gamma_j\circ p_G\circ a_{\Sigma_{\rm\bf{i}}}\\
&=& \gamma_j\circ q\circ(\iota\circ \zeta^{u,e})\circ a_{\Sigma_{\rm\bf{i}}}\\
&=& \gamma_j\circ (\iota\circ \zeta^{u,e})\circ a_{\Sigma_{\rm\bf{i}}}.
\end{eqnarray*}
Thus, we get the first commutative diagram in (\ref{comp2}),
and one can verify other diagrams via similar ways 
for $\Sigma=\Sigma_{\bf{{\rm i}}}$.

Next, we assume that $\Sigma$ is obtained from $\Sigma_{\bf{{\rm i}}}$ by
an iteration of mutations $\ovl{\mu}$. 
Let $\ovl{\mu}^a:\mathcal{A}_{\Sigma_{\bf{{\rm i}}}}\rightarrow \mathcal{A}_{\Sigma}$
(resp. $\ovl{\mu}^x:\mathcal{X}_{\Sigma_{\bf{{\rm i}}}}\rightarrow \mathcal{X}_{\Sigma}$)
be the corresponding birational map.
By Proposition \ref{map-p}, (\ref{commp}) and Definition \ref{actions},
we obtain $a_{\Sigma_{\bf{{\rm i}}}}=a_{\Sigma}\circ\ovl{\mu}^a$,
$x_{\Sigma_{\bf{{\rm i}}}}=x_{\Sigma}\circ\ovl{\mu}^x$ and
$p_{\Sigma_{\bf{{\rm i}}}}=(\ovl{\mu}^x)^{-1}\circ p_{\Sigma}\circ \ovl{\mu}^a$. 
In conjunction with commutative diagrams for $\Sigma_{\bf{{\rm i}}}$,
we obtain the diagrams (\ref{comp1}), (\ref{comp2}) for general seeds $\Sigma$. \qed

\begin{rem}
In the rest of article,
we will treat geometric crystals $(\mathcal{X}_{\Sigma},\ x_{\Sigma}^{-1}\circ e_j^c\circ x_{\Sigma},
\ \varepsilon_j\circ x_{\Sigma},
\ \gamma_j\circ x_{\Sigma})$ and
$(\mathcal{A}_{\Sigma},\ a_{\Sigma}^{-1}\circ(\iota\circ \zeta^{u,e})^{-1}\circ e_j^c
\circ (\iota\circ \zeta^{u,e})\circ a_{\Sigma},
\ \varepsilon_j\circ (\iota\circ \zeta^{u,e})\circ a_{\Sigma},
\ \gamma_j\circ (\iota\circ \zeta^{u,e})\circ a_{\Sigma})$ only.
\end{rem}

\begin{prop}\label{twistprop}
For $i\in [1,\tilde{r}]$ and $w\in W$, we have
\[
(\iota\circ \zeta^{u,e})^{*} \Delta_{w\Lm_i,\Lm_i,}=\Delta_{u\Lm_i,w\Lm_i}.
\]

\end{prop}

\nd
{\it Proof.}

Recall that we defined anti-involution $\omega:G\rightarrow G$ in (\ref{bilin}), One can verify that $\iota\circ \theta=\omega$.
For $x\in G^{u,e}$, we get
\begin{eqnarray*}
\iota\circ\zeta^{u,e}(x)&=&\iota\circ\theta([\ovl{u}^{-1}x]_-^{-1}\ovl{u}^{-1}x)\\
&=&\omega([\ovl{u}^{-1}x]_0[\ovl{u}^{-1}x]_+)\\
&=&\omega([\ovl{u}^{-1}x]_+) \cdot \omega([\ovl{u}^{-1}x]_0).
\end{eqnarray*}

For $x\in G^{u,e}$, we get $\ovl{u}^{-1}x\in G_0$ by Proposition \ref{twistbire}. 
Writing $\ovl{u}^{-1}x=[\ovl{u}^{-1}x]_-[\ovl{u}^{-1}x]_0[\ovl{u}^{-1}x]_+$,
we have $\omega(\ovl{u}^{-1}x)=\omega([\ovl{u}^{-1}x]_+) \omega([\ovl{u}^{-1}x]_0)\omega([\ovl{u}^{-1}x]_-)$
and $\omega([\ovl{u}^{-1}x]_+)\in N^-$, $\omega([\ovl{u}^{-1}x]_0)\in H$, $\omega([\ovl{u}^{-1}x]_-)\in N$.
Using the bilinear form in (\ref{bilin}), we get
\begin{eqnarray*}
(\iota\circ \zeta^{u,e})^{*} \Delta_{w\Lm_i,\Lm_i,}(x)&=&\lan \ovl{w}v_{\Lm_i}, (\iota\circ \zeta^{u,e}) (x) v_{\Lm_i} \ran\\
&=&\lan \ovl{w}v_{\Lm_i}, \omega([\ovl{u}^{-1}x]_+) \cdot \omega([\ovl{u}^{-1}x]_0) v_{\Lm_i} \ran\\
&=&\lan \ovl{w}v_{\Lm_i}, \omega(\ovl{u}^{-1}x) v_{\Lm_i} \ran\\
&=&\lan \ovl{u}^{-1}x \ovl{w}v_{\Lm_i}, v_{\Lm_i} \ran\\
&=&\lan x \ovl{w}v_{\Lm_i}, \ovl{u}v_{\Lm_i}\ran=\Delta_{u\Lm_i,w\Lm_i}(x).
\end{eqnarray*}
\qed

A result similar to Proposition \ref{twistprop} for unipotent quantum minors of quantum unipotent cells is obtained in \cite{KO}.

\section{Explicit formulae of geometric crystals on cluster $\mathcal{X}$-tori}\label{action-X}

In this section, we will reveal the explicit formulae of geometric crystal structures
$(\mathcal{X}_{\Sigma_{\bf{{\rm i}}}},\ x_{\Sigma_{\bf{{\rm i}}}}^{-1}\circ e_j^c\circ x_{\Sigma_{\bf{{\rm i}}}},
\ \varepsilon_j\circ x_{\Sigma_{\bf{{\rm i}}}},
\ \gamma_j\circ x_{\Sigma_{\bf{{\rm i}}}})$ in Definition \ref{actions}.

\begin{thm}\label{thm1}
Let $e^c_j$ $(j\in[1,r],\ c\in\mathbb{C}^{\times})$ be the rational $\mathbb{C}^{\times}$-action on
$G^{u,e}_{{\rm Ad}}$ (Definition \ref{GCdefX}). We set
\[
(X'_{-\tilde{r}},\cdots,X'_{-1},X'_1,\cdots,X'_n):= 
x_{\Sigma_{\bf{{\rm i}}}}^{-1}\circ e^c_j\circ x_{\Sigma_{\bf{{\rm i}}}}(X_{-\tilde{r}},\cdots,X_{-1},X_1,\cdots,X_n),
\]
and $\{K_1,K_2,\cdots,K_l\}:=\{K|1\leq K\leq n ,\ i_K=j\}$ $(K_1<\cdots<K_l)$.
Then
\begin{equation}\label{claim1}
X_{K_p}'
=X_{K_p}
\cdot\frac{c\sum^{p+1}_{m=1}(X_{K_m}X_{K_{m+1}}\cdots X_{K_{l-1}})
+\sum^{l}_{m=p+2}(X_{K_m}X_{K_{m+1}}\cdots X_{K_{l-1}})}{c\sum^{p-1}_{m=1}(X_{K_m}X_{K_{m+1}}\cdots X_{K_{l-1}})+\sum^{l}_{m=p}(X_{K_m}X_{K_{m+1}}\cdots X_{K_{l-1}})}.
\end{equation} 
For $k\in\{1,2,\cdots,n\}$ with $i_k\neq j$, we also set
$\{k_1,k_2,\cdots,k_{s}\}:=\{K|k<K<k^{+},\ j=i_K\}$ $(k_1<k_2<\cdots<k_s)$.
We can write $k_1=K_{\gamma}$ with some $\gamma\in\{1,2,\cdots,l\}$.
Then
\begin{equation}\label{claim2}
X_k'
= X_k\left(\frac{c\sum^{\gamma+s-1}_{m=1}X_{K_{m}}X_{K_{m+1}}\cdots
X_{K_{l-1}}+\sum^{l}_{m=\gamma+s} X_{K_{m}}X_{K_{m+1}}\cdots X_{K_{l-1}}}
{c\sum^{\gamma-1}_{m=1}X_{K_{m}}X_{K_{m+1}}\cdots X_{K_{l-1}}+\sum^{l}_{m=\gamma}X_{K_{m}}X_{K_{m+1}}\cdots X_{K_{l-1}}}
\right)^{a_{j,i_k}}.
\end{equation}
For $i\in[1,\tilde{r}]$,
\[
X_{-i}'=X_{-i}c^{a_{ji}}\prod_{1\leq s\leq n,i_s=i}X_s X_s'^{-1}.
\]
Furthermore, let $\gamma_j$, $\varepsilon_j$ be the functions of
geometric crystal on $G^{u,e}_{{\rm Ad}}$. Then
\[
\gamma_j\circ x_{\Sigma_{\bf{{\rm i}}}}(X_{-\tilde{r}},\cdots,X_{-1},X_1,\cdots,X_n)=X_{-j}X_{K_1}\cdots X_{K_l},
\]
\[
\varepsilon_j\circ x_{\Sigma_{\bf{{\rm i}}}}(X_{-\tilde{r}},\cdots,X_{-1},X_1,\cdots,X_n)
=\left(\sum^{l-1}_{p=0}X_{K_{p+1}}X_{K_{p+2}}\cdots X_{K_l}\right)^{-1}.
\]
\end{thm}

First, let us prove the following lemma.

\begin{lem}\label{sec1-lem}
\begin{enumerate}
\item[(1)] A map $
f:(\mathbb{C}^{\times})^n \rightarrow (\mathbb{C}^{\times})^n$,
\[
(t_1,\cdots,t_s,\cdots,t_n)\mapsto 
(\cdots,\frac{t_{s+1}^{-a_{i_{s+1},i_s}}t_{s+2}^{-a_{i_{s+2},i_s}}\cdots t_{n}^{-a_{i_{n},i_s}}}{t_s},\cdots,\frac{1}{t_n})
\]
is bijective.
\item[(2)]  A map $
g:(\mathbb{C}^{\times})^n \rightarrow (\mathbb{C}^{\times})^n$,
\[
(t_1,\cdots,t_j,\cdots,t_n)\mapsto 
(\prod_{1\leq k\leq n, i_k=i_1}t_k,\cdots,\prod_{j\leq k\leq n, i_k=i_j}t_k, \cdots,t_n)
\]
is bijective.
\end{enumerate}
\end{lem}

\nd
{\it Proof.}

(1) For $(\zeta_1,\cdots,\zeta_n)\in (\mathbb{C}^{\times})^n$,
by setting $t_n=\frac{1}{\zeta_n}$, $t_{n-1}=\frac{(t_n)^{a_{i_n,i_{n-1}}}}{\zeta_{n-1}}$,
$\cdots$, $t_s=\frac{t_{s+1}^{-a_{i_{s+1},i_s}}t_{s+2}^{-a_{i_{s+2},i_s}}\cdots t_{n}^{-a_{i_{n},i_s}}}{\zeta_s}$,
$\cdots$ inductively, we have
\[
(\zeta_1,\cdots,\zeta_n)=f(t_1,\cdots,t_n),
\]
which means $f$ is surjective.

Next, we assume $f(t_1,\cdots,t_n)=f(t_1'\cdots,t_n')$.
Then we get $\frac{1}{t_n}=\frac{1}{t_n'}$, $\frac{t_n^{-a_{i_n,i_{n-1}}}}{t_{n-1}}=\frac{{t'}_n^{-a_{i_n,i_{n-1}}}}{t'_{n-1}}$,
$\cdots$,
$\frac{t_{s+1}^{-a_{i_{s+1},i_s}}t_{s+2}^{-a_{i_{s+2},i_s}}\cdots t_{n}^{-a_{i_{n},i_s}}}{t_s}
=\frac{{t'}_{s+1}^{-a_{i_{s+1},i_s}}{t'}_{s+2}^{-a_{i_{s+2},i_s}}\cdots {t'}_{n}^{-a_{i_{n},i_s}}}{t_s'}$, $\cdots$.
Hence, we can inductively show that $t_n=t_n'$, $t_{n-1}=t_{n-1}'$, $\cdots$, $t_s=t_s'$, $\cdots$.
Therefore, $f$ is injective. Similarly, we can prove (2). \qed

\vspace{3mm}

\nd
{\it Proof of Theorem \ref{thm1}.} 

Note that
$y_i(t)X^{\Lm^{\vee}_l}=X^{\Lm^{\vee}_l}y_i(X^{\delta_{i,l}} t)$ holds.
Therefore,
\begin{eqnarray*}
& &x_{\Sigma_{\bf{{\rm i}}}}(X_{-\tilde{r}},\cdots,X_{-1},X_1,\cdots,X_n)\\
&=&X_{-\tilde{r}}^{\Lm_{\tilde{r}}^{\vee}}\cdots X_{-1}^{\Lm_{1}^{\vee}}
y_{i_1}(1) X_1^{\Lm_{i_1}^{\vee}} y_{i_2}(1) X_2^{\Lm_{i_2}^{\vee}}\cdots y_{i_n}(1) X_n^{\Lm_{i_n}^{\vee}}\\
&=& X_{-\tilde{r}}^{\Lm_{\tilde{r}}^{\vee}}\cdots X_{-1}^{\Lm_{1}^{\vee}}  X_1^{\Lm_{i_1}^{\vee}} \cdots X_n^{\Lm_{i_n}^{\vee}}\\
& &\cdot y_{i_1}(\prod_{1\leq k\leq n,i_k=i_1}X_k) y_{i_2}(\prod_{2\leq k\leq n,i_k=i_2}X_k)
y_{i_3}(\prod_{3\leq k\leq n,i_k=i_3}X_k)\cdots
y_{i_n}(X_n).
\end{eqnarray*}

By Lemma \ref{sec1-lem},
there exist $t_1,\cdots,t_n\in\mathbb{C}^{\times}$ such that
\begin{equation}\label{1-00}
\prod_{s\leq k\leq n,i_k=i_s}X_k =\frac{t_{s+1}^{-a_{i_{s+1},i_s}}t_{s+2}^{-a_{i_{s+2},i_s}}\cdots t_{n}^{-a_{i_{n},i_s}}}{t_s}
\end{equation}
for $s=1,2,\cdots,n$.

Hence,
\begin{eqnarray}
& &x_{\Sigma_{\bf{{\rm i}}}}(X_{-\tilde{r}},\cdots,X_{-1},X_1,\cdots,X_n)\nonumber\\
&=& X_{-\tilde{r}}^{\Lm_{\tilde{r}}^{\vee}}\cdots X_{-1}^{\Lm_{1}^{\vee}}  X_1^{\Lm_{i_1}^{\vee}} \cdots X_n^{\Lm_{i_n}^{\vee}}\nonumber\\
& &\cdot y_{i_1}(\frac{t_2^{-a_{i_2,i_1}}t_3^{-a_{i_3,i_1}}\cdots t_n^{-a_{i_n,i_1}}}{t_1}) 
y_{i_2}(\frac{t_3^{-a_{i_3,i_2}}t_4^{-a_{i_4,i_2}}\cdots t_n^{-a_{i_n,i_2}}}{t_2})\nonumber\\
& &\cdot y_{i_3}(\frac{t_4^{-a_{i_4,i_3}}\cdots t_n^{-a_{i_n,i_3}}}{t_3})\cdots
y_{i_n}(\frac{1}{t_n})\nonumber\\
&=& X_{-\tilde{r}}^{\Lm_{\tilde{r}}^{\vee}}\cdots X_{-1}^{\Lm_{1}^{\vee}}  X_1^{\Lm_{i_1}^{\vee}} \cdots X_n^{\Lm_{i_n}^{\vee}}
\al^{\vee}_{i_1}(t_1)\al^{\vee}_{i_2}(t_2)\cdots \al^{\vee}_{i_n}(t_n)\label{1-0}
\\
& &\cdot y_{i_1}(t_1)\al^{\vee}_{i_1}(t_1^{-1}) y_{i_2}(t_2)\al^{\vee}_{i_2}(t_2^{-1})\cdots y_{i_n}(t_n)\al^{\vee}_{i_n}(t_n^{-1})\nonumber,
\end{eqnarray}
where we use
\begin{equation}\label{1-1}
y_{i}(S)\al^{\vee}_{p}(T^{-1})=\al^{\vee}_p(T^{-1})y_i(S T^{-a_{p,i}}),\quad 
(S,T\in\mathbb{C}^{\times},\ i,p\in\{1,2,\cdots,r\})
\end{equation}
 in the second equality.
Here, we denote the image of $\al^{\vee}_p(T^{-1})\in G^{u,e}$ under the quotient map $G\rightarrow G_{\rm Ad}$
 by the same notation $\al^{\vee}_p(T^{-1})$.

By the definition of $e^c_j$ and Proposition \ref{geom-prop},
if $(c,x_{\Sigma_{\bf{{\rm i}}}}(X_{-\tilde{r}},\cdots,X_{-1},X_1,\cdots,X_n))\in {\rm dom}(e_j)$ then we have
\begin{eqnarray*}
& &
e^c_j \circ x_{\Sigma_{\bf{{\rm i}}}}(X_{-\tilde{r}},\cdots,X_{-1},X_1,\cdots,X_n)\\
&=& X_{-\tilde{r}}^{\Lm_{\tilde{r}}^{\vee}}\cdots X_{-1}^{\Lm_{1}^{\vee}}  X_1^{\Lm_{i_1}^{\vee}} \cdots X_n^{\Lm_{i_n}^{\vee}}
\al^{\vee}_{i_1}(t_1)\al^{\vee}_{i_2}(t_2)\cdots \al^{\vee}_{i_n}(t_n)
\\
& &\cdot y_{i_1}(t_1')\al^{\vee}_{i_1}(t_1'^{-1}) y_{i_2}(t_2')\al^{\vee}_{i_2}(t_2'^{-1})\cdots y_{i_n}(t_n')\al^{\vee}_{i_n}(t_n'^{-1}),
\end{eqnarray*}
with $t_1',\cdots,t_n'\in\mathbb{C}^{\times}$ in Proposition \ref{geom-prop}.

Using (\ref{1-1}) again, we obtain
\begin{eqnarray}
& &
e^c_j \circ x_{\Sigma_{\bf{{\rm i}}}}(X_{-\tilde{r}},\cdots,X_{-1},X_1,\cdots,X_n)\nonumber\\
&=& X_{-\tilde{r}}^{\Lm_{\tilde{r}}^{\vee}}\cdots X_{-1}^{\Lm_{1}^{\vee}}  X_1^{\Lm_{i_1}^{\vee}} \cdots X_n^{\Lm_{i_n}^{\vee}}
\al^{\vee}_{i_1}(t_1)\al^{\vee}_{i_2}(t_2)\cdots \al^{\vee}_{i_n}(t_n)
\nonumber\\
& &\cdot \al^{\vee}_{i_1}(t_1'^{-1}) \al^{\vee}_{i_2}(t_2'^{-1})\cdots \al^{\vee}_{i_n}(t_n'^{-1})
y_{i_1}(\frac{{t'}_2^{-a_{i_2,i_1}}{t'}_3^{-a_{i_3,i_1}}\cdots {t'}_n^{-a_{i_n,i_1}}}{t_1'}) \label{eq2}\\
& &
\cdot y_{i_2}(\frac{{t'}_3^{-a_{i_3,i_2}}{t'}_4^{-a_{i_4,i_2}}\cdots {t'}_n^{-a_{i_n,i_2}}}{t_2'})
 y_{i_3}(\frac{{t'}_4^{-a_{i_4,i_3}}\cdots {t'}_n^{-a_{i_n,i_3}}}{t_3'})\cdots
y_{i_n}(\frac{1}{t_n'}).\nonumber
\end{eqnarray}

It follows from Lemma \ref{sec1-lem} that
there exist $X_1',\cdots,X_n'\in\mathbb{C}^{\times}$ such that
\begin{equation}\label{eq3}
\frac{(t_{k+1}')^{-a_{i_{k+1},i_k}}(t_{k+2}')^{-a_{i_{k+2},i_k}}\cdots (t_n')^{-a_{i_{n},i_k}}}{t_k'}
=\prod_{k\leq p\leq n,i_p=i_k}X_p'
\end{equation}
for $k=1,2,\cdots,n$.
Substituting these into (\ref{eq2}), we get
\begin{eqnarray*}
& &
e^c_j \circ x_{\Sigma_{\bf{{\rm i}}}}(X_{-\tilde{r}},\cdots,X_{-1},X_1,\cdots,X_n)\\
&=& X_{-\tilde{r}}^{\Lm_{\tilde{r}}^{\vee}}\cdots X_{-1}^{\Lm_{1}^{\vee}}  X_1^{\Lm_{i_1}^{\vee}} \cdots X_n^{\Lm_{i_n}^{\vee}}
\al^{\vee}_{i_1}(t_1)\al^{\vee}_{i_2}(t_2)\cdots \al^{\vee}_{i_n}(t_n)
\\
& &\cdot \al^{\vee}_{i_1}(t_1'^{-1}) \al^{\vee}_{i_2}(t_2'^{-1})\cdots \al^{\vee}_{i_n}(t_n'^{-1})
y_{i_1}(\prod_{1\leq l\leq n,i_l=i_1}X_l') \\
& &
\cdot y_{i_2}(\prod_{2\leq l\leq n,i_l=i_2}X_l')
 y_{i_3}(\prod_{3\leq l\leq n,i_l=i_3}X_l')\cdots
y_{i_n}(X_n')\\
&=& X_{-\tilde{r}}^{\Lm_{\tilde{r}}^{\vee}}\cdots X_{-1}^{\Lm_{1}^{\vee}}  X_1^{\Lm_{i_1}^{\vee}} \cdots X_n^{\Lm_{i_n}^{\vee}}
\al^{\vee}_{i_1}(t_1)\al^{\vee}_{i_2}(t_2)\cdots \al^{\vee}_{i_n}(t_n)
\\
& &\cdot \al^{\vee}_{i_1}(t_1'^{-1}) \al^{\vee}_{i_2}(t_2'^{-1})\cdots \al^{\vee}_{i_n}(t_n'^{-1})
 X_1'^{-\Lm_{i_1}^{\vee}} \cdots X_n'^{-\Lm_{i_n}^{\vee}}\\
& &y_{i_1}(1)X_1'^{\Lm_{i_1}^{\vee}}  
y_{i_2}(1) X_2'^{\Lm_{i_2}^{\vee}} \cdots
y_{i_n}(1)X_n'^{\Lm_{i_n}^{\vee}},
\end{eqnarray*}
where we use $y_i(T)X^{\Lm^{\vee}_l}=X^{\Lm^{\vee}_l}y_i(X^{\delta_{i,l}} T)$ in the second equality.
Note that putting $\{K_1,K_2,\cdots,K_l\}:=\{1\leq K\leq n| i_K=j\}$\ ($K_1<\cdots<K_l$),
\begin{eqnarray}
& & t_k t_k'^{-1}=\nonumber\\ 
& &\hspace{-8mm}\begin{cases}
\frac{\displaystyle c\sum^{p}_{m=1}t_1^{a_{i_1,j}}t_2^{a_{i_2,j}}\cdots t_{K_m-1}^{a_{i_{K_m-1},j}}t_{K_m}+\sum^{l}_{m=p+1}t_1^{a_{i_1,j}}t_2^{a_{i_2,j}}\cdots  t_{K_m-1}^{a_{i_{K_m-1},j}}t_{K_m}}{\displaystyle c\sum^{p-1}_{m=1}t_1^{a_{i_1,j}}t_2^{a_{i_2,j}}\cdots  t_{K_m-1}^{a_{i_{K_m-1},j}}t_{K_m}+\sum^{l}_{m=p}t_1^{a_{i_1,j}}t_2^{a_{i_2,j}}\cdots  t_{K_m-1}^{a_{i_{K_m-1},j}}t_{K_m}} & {\rm if}\ k=K_p,
 \label{eq3-1} \\
1 & {\rm if}\ i_k\neq j.
\end{cases}
\end{eqnarray}

Thus, by Proposition \ref{ExDyn}, we get
\begin{eqnarray*}
& &\al^{\vee}_{i_1}(t_1)\al^{\vee}_{i_2}(t_2)\cdots \al^{\vee}_{i_n}(t_n)
\cdot \al^{\vee}_{i_1}(t_1'^{-1}) \al^{\vee}_{i_2}(t_2'^{-1})\cdots \al^{\vee}_{i_n}(t_n'^{-1})\\
&=& \al^{\vee}_j (t_{K_1}t_{K_1}'^{-1}t_{K_2}t_{K_2}'^{-1}\cdots t_{K_l}t_{K_l}'^{-1})=\al^{\vee}_j(c)=c^{\sum^{\tilde{r}}_{i=1}a_{ji}\Lm^{\vee}_i}.
\end{eqnarray*}

Therefore, it follows that
\begin{eqnarray*}
& &
e^c_j \circ x_{\Sigma_{\bf{{\rm i}}}}(X_{-\tilde{r}},\cdots,X_{-1},X_1,\cdots,X_n)\\
&=& X_{-\tilde{r}}'^{\Lm_{\tilde{r}}^{\vee}}\cdots X_{-1}'^{\Lm_{1}^{\vee}}  
 y_{i_1}(1)X_1'^{\Lm_{i_1}^{\vee}}  
y_{i_2}(1) X_2'^{\Lm_{i_2}^{\vee}} \cdots
y_{i_n}(1)X_n'^{\Lm_{i_n}^{\vee}},
\end{eqnarray*}
where for $i\in[1,\tilde{r}]$,
\[
X_{-i}'=X_{-i}c^{a_{ji}}\prod_{1\leq s\leq n,i_s=i}X_s X_s'^{-1}.
\]
Finally, let us prove (\ref{claim1}) and (\ref{claim2}).
For $p=1,\cdots,l$, the definition $(\ref{eq3})$ implies that
\begin{equation}\label{eq4}
X_{K_p}'=
 \frac{(t_{K_p+1}')^{-a_{i_{K_p+1},j}}(t_{K_p+2}')^{-a_{i_{K_p+2},j}}\cdots (t_n')^{-a_{i_{n},j}}}{X_{K_{p+1}}'\cdots X_{K_l}'t_{K_p}'}.
\end{equation}

By using induction on $l-p$, let us prove
\begin{equation}\label{eq5}
X_{K_p}'=\frac{{t'}_{K_p+1}^{-a_{i_{K_p+1},j}}{t'}_{K_p+2}^{-a_{i_{K_p+2},j}}\cdots {t'}_{K_{p+1}-1}^{-a_{i_{K_{p+1}-1},j}}}{t_{K_p}'t_{K_{p+1}}'},\ \ (p=1,2,\cdots,l)
\end{equation}
where $T'_{K_{l+1}}:=1$.
First, we consider the case $l-p=0$ so that $p=l$.
In this case, our claim follows from (\ref{eq4}).
Next, we assume that $l>p$.
Taking into account the hypothesis of induction, by $t_{K_{p+1}}'t_{K_{p+2}}'^2
t_{K_{p+3}}'^2\cdots t_{K_{l}}'^2=\frac{1}{{t'}_{K_{p+1}}}{t'}_{K_{p+1}}^{a_{i_{K_{p+1}},j}}{t'}_{K_{p+2}}^{a_{i_{K_{p+2}},j}}{t'}_{K_{p+3}}^{a_{i_{K_{p+3}},j}}\cdots {t'}_{K_{l}}^{a_{i_{K_{l}},j}}$,
\begin{eqnarray*}
& &X_{K_{p+1}}'\cdots X_{K_l}'= \frac{1}{t_{K_{p+1}}'t_{K_{p+2}}'^2
t_{K_{p+3}}'^2\cdots t_{K_{l}}'^2} {t'}_{K_{p+1}+1}^{-a_{i_{K_{p+1}+1},j}}{t'}_{K_{p+1}+2}^{-a_{i_{K_{p+1}+2},j}}\cdots {t'}_{K_{p+2}-1}^{-a_{i_{K_{p+2}-1},j}}\\
& &
{t'}_{K_{p+2}+1}^{-a_{i_{K_{p+2}+1},j}}{t'}_{K_{p+2}+2}^{-a_{i_{K_{p+2}+2},j}}\cdots {t'}_{K_{p+3}-1}^{-a_{i_{K_{p+3}-1},j}}\cdots {t'}_{K_l+1}^{-a_{i_{K_l+1},j}}{t'}_{K_l+2}^{-a_{i_{K_l+2},j}}\cdots {t'}_{n}^{-a_{i_{n},j}}\\
&=& t_{K_{p+1}}' {t'}_{K_{p+1}}^{-a_{i_{K_{p+1}},j}}{t'}_{K_{p+1}+1}^{-a_{i_{K_{p+1}+1},j}}{t'}_{K_{p+1}+2}^{-a_{i_{K_{p+1}+2},j}}\cdots {t'}_{n}^{-a_{i_{n},j}}.
\end{eqnarray*}
In conjunction with (\ref{eq4}), we obtain (\ref{eq5}).
Using (\ref{eq3-1}), we obtain
\begin{eqnarray}
X_{K_p}'
&=&\frac{t_{K_p+1}^{-a_{i_{K_p+1},j}}t_{K_p+2}^{-a_{i_{K_p+2},j}}\cdots t_{K_{p+1}-1}^{-a_{i_{K_{p+1}-1},j}}}{t_{K_p}'t_{K_{p+1}}'}\nonumber\\
&=&\frac{t_{K_p+1}^{-a_{i_{K_p+1},j}}t_{K_p+2}^{-a_{i_{K_p+2},j}}\cdots t_{K_{p+1}-1}^{-a_{i_{K_{p+1}-1},j}}}{t_{K_p}t_{K_{p+1}}}\label{eq6}\\
& &
\cdot\frac{c\sum^{p+1}_{m=1}t_1^{a_{i_1,j}}t_2^{a_{i_2,j}}\cdots t_{K_m-1}^{a_{i_{K_m-1},j}}t_{K_m}+\sum^{l}_{m=p+2}t_1^{a_{i_1,j}}t_2^{a_{i_2,j}}\cdots  t_{K_m-1}^{a_{i_{K_m-1},j}}t_{K_m}}{c\sum^{p-1}_{m=1}t_1^{a_{i_1,j}}t_2^{a_{i_2,j}}\cdots  t_{K_m-1}^{a_{i_{K_m-1},j}}t_{K_m}+\sum^{l}_{m=p}t_1^{a_{i_1,j}}t_2^{a_{i_2,j}}\cdots  t_{K_m-1}^{a_{i_{K_m-1},j}}t_{K_m}}.\nonumber
\end{eqnarray}
Substituting $c=1$, we obtain
\begin{equation}\label{eq7}
X_{K_p}=
\frac{t_{K_p+1}^{-a_{i_{K_p+1},j}}t_{K_p+2}^{-a_{i_{K_p+2},j}}\cdots t_{K_{p+1}-1}^{-a_{i_{K_{p+1}-1},j}}}{t_{K_p}t_{K_{p+1}}}.
\end{equation}
Therefore, 
\[
(X_{K_1}X_{K_2}\cdots X_{K_{m-1}})^{-1}=
t_{K_1}t_{K_1+1}^{a_{i_{K_1+1},j}}t_{K_1+2}^{a_{i_{K_1+2},j}}t_{K_1+3}^{a_{i_{K_1+3},j}}\cdots t_{K_m-1}^{a_{i_{K_m-1},j}}t_{K_m},
\]
where we use $t_{K_p}^{a_{i_{K_p},j}}=t_{K_p}^{a_{j,j}}=t_{K_p}^2$. Thus,
\[
t_1^{a_{i_1,j}}t_2^{a_{i_2,j}}\cdots  t_{K_1-1}^{a_{i_{K_1-1},j}}t_{K_1}
(X_{K_1}X_{K_2}\cdots X_{K_{m-1}})^{-1}
=t_1^{a_{i_1,j}}t_2^{a_{i_2,j}}t_3^{a_{i_3,j}}\cdots t_{K_m-1}^{a_{i_{K_m-1},j}}t_{K_m}.
\]
It follows from (\ref{eq6}) that
\begin{eqnarray*}
& &X_{K_p}'\\
&=&X_{K_p}
\cdot\frac{c\sum^{p+1}_{m=1}(X_{K_1}X_{K_2}\cdots X_{K_{m-1}})^{-1}
+\sum^{l}_{m=p+2}(X_{K_1}X_{K_2}\cdots X_{K_{m-1}})^{-1}}{c\sum^{p-1}_{m=1}(X_{K_1}X_{K_2}\cdots X_{K_{m-1}})^{-1}+\sum^{l}_{m=p}(X_{K_1}X_{K_2}\cdots X_{K_{m-1}})^{-1}}\\
&=&X_{K_p}
\cdot\frac{c\sum^{p+1}_{m=1}(X_{K_m}X_{K_{m+1}}\cdots X_{K_{l-1}})
+\sum^{l}_{m=p+2}(X_{K_m}X_{K_{m+1}}\cdots X_{K_{l-1}})}{c\sum^{p-1}_{m=1}(X_{K_m}X_{K_{m+1}}\cdots X_{K_{l-1}})+\sum^{l}_{m=p}(X_{K_m}X_{K_{m+1}}\cdots X_{K_{l-1}})}.
\end{eqnarray*}
Thus, we obtain (\ref{claim1}).

Next, for $k\in\{1,2,\cdots,n\}$, let us suppose that $i_k\neq j$.
We set $\{k_1,k_2,\cdots,k_{s}\}:=\{K|k<K<k^{+},\ j=i_K\}$ $(k_1<k_2<\cdots<k_s)$.
Because of (\ref{eq3}), 
\begin{eqnarray*}
& &X_k'=\frac{\prod_{k\leq p\leq n,i_p=i_k}X_p'}{\prod_{k^+\leq p\leq n,i_p=i_k}X_p'}\\
&=& \frac{(t_{k+1}')^{-a_{i_{k+1},i_k}}(t_{k+2}')^{-a_{i_{k+2},i_k}}\cdots (t_{k^+ -1}')^{-a_{i_{k^+ -1},i_k}}}{t_k't_{k^+}'}\\
&=& \frac{1}{t_kt_{k^+}}(t_{k+1})^{-a_{i_{k+1},i_k}}\cdots (t_{k_1-1})^{-a_{i_{k_1-1},i_k}} (t_{k_1}')^{-a_{j,i_k}}
 (t_{k_1+1})^{-a_{i_{k_1+1},i_k}} \\
& &\cdots  (t_{k_2-1})^{-a_{i_{k_2-1},i_k}} (t_{k_2}')^{-a_{j,i_k}}
 (t_{k_2+1})^{-a_{i_{k_2+1},i_k}}\\
& &\cdots  (t_{k_s-1})^{-a_{i_{k_s-1},i_k}} (t_{k_s}')^{-a_{j,i_k}}
 (t_{k_s+1})^{-a_{i_{k_s+1},i_k}}
\cdots (t_{k^+ -1})^{-a_{i_{k^+ -1},i_k}}.
\end{eqnarray*}

Since $i_{k_1}=j$, there exists $\gamma\in\{1,2,\cdots,l\}$ such that $k_1=K_{\gamma}$.
In this case, we have $k_s=K_{\gamma+s-1}$. By (\ref{eq3-1}) and (\ref{eq7}),
\begin{eqnarray*}
& &X_k' \\
&=&X_k\left(\frac{c\sum^{\gamma+s-1}_{m=1}t_1^{a_{i_1,j}}t_2^{a_{i_2,j}}\cdots t_{K_m-1}^{a_{i_{K_m-1},j}}t_{K_m}+\sum^{l}_{m=\gamma+s}t_1^{a_{i_1,j}}t_2^{a_{i_2,j}}\cdots  t_{K_m-1}^{a_{i_{K_m-1},j}}t_{K_m}}{c\sum^{\gamma-1}_{m=1}t_1^{a_{i_1,j}}t_2^{a_{i_2,j}}\cdots  t_{K_m-1}^{a_{i_{K_m-1},j}}t_{K_m}+\sum^{l}_{m=\gamma}t_1^{a_{i_1,j}}t_2^{a_{i_2,j}}\cdots  t_{K_m-1}^{a_{i_{K_m-1},j}}t_{K_m}}\right)^{a_{j,i_k}}\\
&=& X_k\left(\frac{c\sum^{\gamma+s-1}_{m=1}(X_{K_1}X_{K_2}\cdots X_{K_{m-1}})^{-1}+\sum^{l}_{m=\gamma+s}(X_{K_1}X_{K_2}\cdots X_{K_{m-1}})^{-1}}{c\sum^{\gamma-1}_{m=1}(X_{K_1}X_{K_2}\cdots X_{K_{m-1}})^{-1}+\sum^{l}_{m=\gamma}(X_{K_1}X_{K_2}\cdots X_{K_{m-1}})^{-1}}
\right)^{a_{j,i_k}}\\
&=& X_k\left(\frac{c\sum^{\gamma+s-1}_{m=1}X_{K_{m}}X_{K_{m+1}}\cdots X_{K_{l-1}}+\sum^{l}_{m=\gamma+s} X_{K_{m}}X_{K_{m+1}}\cdots X_{K_{l-1}}}{c\sum^{\gamma-1}_{m=1}X_{K_{m}}X_{K_{m+1}}\cdots X_{K_{l-1}}+\sum^{l}_{m=\gamma}X_{K_{m}}X_{K_{m+1}}\cdots X_{K_{l-1}}}
\right)^{a_{j,i_k}}.
\end{eqnarray*}
Thus, we get (\ref{claim2}).

Taking into account (\ref{1-0}), we get 
\begin{eqnarray*}
\gamma_j\circ x_{\Sigma_{\bf{{\rm i}}}}(X_{-\tilde{r}},\cdots,X_{-1},X_1,\cdots,X_n)&=&
\frac{\al_j( X_{-\tilde{r}}^{\Lm_{\tilde{r}}^{\vee}}\cdots X_{-1}^{\Lm_{1}^{\vee}}  X_1^{\Lm_{i_1}^{\vee}} \cdots X_n^{\Lm_{i_n}^{\vee}}
\al^{\vee}_{i_1}(t_1)\al^{\vee}_{i_2}(t_2)\cdots \al^{\vee}_{i_n}(t_n))}{t_1^{a_{i_1,j}}\cdots t_n^{a_{i_n,j}}}\\
&=&\al_j( X_{-\tilde{r}}^{\Lm_{\tilde{r}}^{\vee}}\cdots X_{-1}^{\Lm_{1}^{\vee}}  X_1^{\Lm_{i_1}^{\vee}} \cdots X_n^{\Lm_{i_n}^{\vee}})\\
&=&X_{-j}X_{K_1}\cdots X_{K_l}.
\end{eqnarray*}
Using (\ref{1-00}), we also get
\begin{eqnarray*}
\varepsilon_j\circ x_{\Sigma_{\bf{{\rm i}}}}(X_{-\tilde{r}},\cdots,X_{-1},X_1,\cdots,X_n)&=&
\left(\sum_{1\leq m\leq n,\ i_m=j}\frac{1}{t_mt_{m+1}^{a_{i_{m+1},j}}\cdots t_{n}^{a_{i_{n},j}}}\right)^{-1}\\
&=&\left(\sum^{l-1}_{p=0} X_{K_{p+1}}X_{K_{p+2}}\cdots X_{K_l} \right)^{-1}.
\end{eqnarray*}

\qed

By Theorem \ref{thm1}, we can verify that the map
$x_{\Sigma_{\rm\bf{i}}}:\mathcal{X}_{\Sigma_{\rm \bf{i}}}\rightarrow G^{u,e}_{{\rm Ad}}$
is a positive structure on the geometric crystal
$(G^{u,e}_{\rm Ad}, \{e_i\}_{i\in [1,r]},\{\gamma_i\}_{i\in [1,r]},\{\varepsilon_i\}_{i\in [1,r]}))$.
Let $\Sigma$ be a seed obtained from $\Sigma_{\rm\bf{i}}$ by an iteration
of mutations. Then the corresponding birational map 
$\mathcal{X}_{\Sigma_{\rm \bf{i}}}\rightarrow \mathcal{X}_{\Sigma}$
is positive. Hence, we obtain the following:

\begin{cor}\label{corX}
The map $x_{\Sigma}:\mathcal{X}_{\Sigma}\rightarrow G^{u,e}_{{\rm Ad}}$ is a
positive structure on the geometric crystal
$(G^{u,e}_{\rm Ad}, \{e_i\}_{i\in [1,r]},\{\gamma_i\}_{i\in [1,r]},\{\varepsilon_i\}_{i\in [1,r]}))$.
Applying
Theorem \ref{UDthm}, we obtain a crystal 
\[(X_{*}(\mathcal{X}_{\Sigma}), \{\tilde{e}_i\}_{i\in [1,r]},\{\tilde{\gamma}_i\}_{i\in [1,r]},\{\tilde{\varepsilon}_i\}_{i\in [1,r]}). \]
It induces a glued crystal structure on the set of $\mathbb{Z}$-valued points $\mathcal{X}_{|\Sigma_{\rm \bf{i}}|}(\mathbb{Z}^T)$.
\end{cor}

\section{Explicit formulae of geometric crystals on cluster $\mathcal{A}$-tori}\label{action-A}

In this section, we will present the explicit formulae of geometric crystal structures
$(\mathcal{A}_{\Sigma_{\bf{{\rm i}}}},\ a_{\Sigma_{\bf{{\rm i}}}}^{-1}\circ(\iota\circ \zeta^{u,e})^{-1}\circ e_j^c
\circ (\iota\circ \zeta^{u,e})\circ a_{\Sigma_{\bf{{\rm i}}}},
\ \varepsilon_j\circ (\iota\circ \zeta^{u,e})\circ a_{\Sigma_{\bf{{\rm i}}}},
\ \gamma_j\circ (\iota\circ \zeta^{u,e})\circ a_{\Sigma_{\bf{{\rm i}}}})$ in Definition \ref{actions}.
For a fixed reduced word $\textbf{i}=(i_1,\cdots,i_n)$ of $u$ and $k\in\{1,2,\cdots,n\}$,
we set $k^-:={\rm max}\{l\in I | l<k, |i_l|=|i_k|\}$.
For $j\in[1,r]$, we also set $j_{\rm max}:={\rm max}\{l\in I| i_l=j\}$.


\begin{thm}\label{thm4}
Let $M$ be the matrix in Theorem \ref{pgdef}, and $p=p_M$ be the map in Proposition \ref{map-p}.
We put
\[
(\overline{X}_{-\tilde{r}},\cdots,\overline{X}_{-1},\overline{X}_{1},\cdots,\overline{X}_{n})
=x_{\Sigma_{\bf{{\rm i}}}}^{-1} \circ e_j^c \circ x_{\Sigma_{\bf{{\rm i}}}} \circ 
 p(A_{-\tilde{r}},\cdots,A_{-1},A_1,\cdots,A_n).
\]
Then for $j\in\{1,\cdots,r\}$ and $k\in \{1,\cdots,n\}$, we have
\begin{eqnarray}
& &(a_{\Sigma_{\bf{{\rm i}}}}^{-1}\circ (\iota\circ\zeta^{u,e})^{-1} \circ e_j^c \circ  (\iota\circ\zeta^{u,e}) \circ a_{\Sigma_{\bf{{\rm i}}}}A)_{k^-} \label{claim5-1}\\
&=& A_{k^-}\cdot c^{D_{\{k,k+1,\cdots,n\},\{(k+1)^-,(k+2)^-,\cdots,n^-,j_{\rm max}\}}}
\left(\frac{p(A)_{k}}{\overline{X}_{k}}\right)
\left(\frac{p(A)_{k+1}}{\overline{X}_{k+1}}\right)^{D_{k,(k+1)^-}} \nonumber\\
& &\cdot 
\left( \frac{p(A)_{k+2}}{\overline{X}_{k+2}}\right)^{D_{\{k,k+1\},\{(k+1)^-,(k+2)^-\}}}
\cdots
\left(\frac{p(A)_{n}}{\overline{X}_{n}}\right)^{D_{\{k,k+1,\cdots,n-1\},\{(k+1)^-,(k+2)^-,\cdots,n^-\}}},\ \nonumber
\end{eqnarray}
where $D_{\{j_1,\cdots,j_l\},\{k_1,\cdots,k_l\}}$ is the minor of $\tilde{B}_{{\rm \bf{i}}}=B_{{\rm \bf{i}}}+M=(\tilde{B}_{i,s})_{i,s\in\{-\tilde{r},\cdots,-1,1,2,\cdots,n\}}$ whose rows (resp. columns) are labelled by
$\{j_1,\cdots,j_l\}$ (resp $\{k_1,\cdots,k_l\}$).
Furthermore, for $k\in \{1,2,\cdots,r\}$, we obtain 
\begin{equation}\label{claim5-2}
(a_{\Sigma_{\bf{{\rm i}}}}^{-1}\circ (\iota\circ\zeta^{u,e})^{-1} \circ e_j^c \circ  (\iota\circ\zeta^{u,e}) \circ a_{\Sigma_{\bf{{\rm i}}}}A)_{k_{\rm max}}=c^{\delta_{j,k}}A_{k_{\rm max}}.
\end{equation}
If $i\in \{1,2,\cdots,\tilde{r}\}\setminus \{i_1,i_2,\cdots,i_n\}$ then
\begin{equation}\label{claim5-3}
(a_{\Sigma_{\bf{{\rm i}}}}^{-1}\circ (\iota\circ\zeta^{u,e})^{-1} \circ e_j^c \circ  (\iota\circ\zeta^{u,e}) \circ a_{\Sigma_{\bf{{\rm i}}}}A)_{-i}=A_{-i}.
\end{equation}
\end{thm}

\begin{rem}
In the above theorem, we obtain formulae
of $(a_{\Sigma_{\bf{{\rm i}}}}^{-1}\circ (\iota\circ\zeta^{u,e})^{-1} \circ 
e_j^c \circ  (\iota\circ\zeta^{u,e}) \circ a_{\Sigma_{\bf{{\rm i}}}}A)_{k}$
for all $k\in I$.
\end{rem}

\nd
{\it Proof of Theorem \ref{thm4}.}

First, let us prove (\ref{claim5-2}).
By a property of the irreducible highest weight module $V(\Lm_k)$ with the highest weight vector
$v_{\Lm_k}$, we obtain $f_j v_{\Lm_k}=0$ $(j\neq k)$ and $f_k v_{\Lm_k}=\ovl{s_k}v_{\Lm_k}$ (\cite{Kum}, Chap.II).
If $j\neq k$, by using the bilinear form in \ref{bilin}, we obtain
\begin{eqnarray*}
((e_j^c)^* \Delta_{\Lm_k,\Lm_k})(x)&=&\Delta_{\Lm_k,\Lm_k}(x_j((c-1)\varphi_j(x)) x x_j((c^{-1}-1)\varepsilon_j(x))) \\
&=&\lan v_{\Lm_k}, x_j((c-1)\varphi_j(x)) x x_j((c^{-1}-1)\varepsilon_j(x)) v_{\Lm_k} \ran\\
&=&\lan y_j((c-1)\varphi_j(x))v_{\Lm_k}, x v_{\Lm_k} \ran\\
&=&\lan v_{\Lm_k}, x v_{\Lm_k} \ran\\
&=&\Delta_{\Lm_k,\Lm_k}(x),
\end{eqnarray*}
 which means $(e_j^c)^* \Delta_{\Lm_k,\Lm_k}=\Delta_{\Lm_k,\Lm_k}$. 

We also obtain
\begin{eqnarray*}
((e_j^c)^* \Delta_{\Lm_j,\Lm_j})(x)
&=&\lan y_j((c-1)\varphi_j(x))v_{\Lm_j}, x v_{\Lm_j} \ran\\
&=&\lan v_{\Lm_j}, x v_{\Lm_j} \ran
+ (c-1)\varphi_j(x) \lan f_jv_{\Lm_j}, x v_{\Lm_j} \ran 
\\
&=&\lan v_{\Lm_j}, x v_{\Lm_j} \ran
+ (c-1)\varphi_j(x) \lan \ovl{s}_jv_{\Lm_j}, x v_{\Lm_j} \ran 
\\
&=&\lan v_{\Lm_j}, x v_{\Lm_j} \ran
+ (c-1)\lan v_{\Lm_j}, x v_{\Lm_j} \ran 
= c\lan v_{\Lm_j}, x v_{\Lm_j} \ran,
\end{eqnarray*}
where we use (\ref{minorexp}) in the fourth equality.
Thus, we get $(e_j^c)^* \Delta_{\Lm_j,\Lm_j}=c\Delta_{\Lm_j,\Lm_j}$ and
\begin{eqnarray*}
& &(a_{\Sigma_{\bf{{\rm i}}}}^{-1}\circ (\iota\circ \zeta^{u,e})^{-1} \circ e_j^c\circ (\iota\circ \zeta^{u,e}) \circ a_{\Sigma_{\bf{{\rm i}}}})^*A_{k_{\rm max}}\\
&=& a_{\Sigma_{\bf{{\rm i}}}}^* \circ (\iota\circ \zeta^{u,e})^{*} \circ (e_j^c)^*\circ
((\iota\circ \zeta^{u,e})^{-1})^*\circ(a_{\Sigma_{\bf{{\rm i}}}}^{-1})^*A_{k_{\rm max}}\\
&=& a_{\Sigma_{\bf{{\rm i}}}}^* \circ (\iota\circ \zeta^{u,e})^{*} \circ (e_j^c)^*\circ
((\iota\circ \zeta^{u,e})^{-1})^* \Delta_{u\Lm_k,\Lm_k} \\
&=& a_{\Sigma_{\bf{{\rm i}}}}^* \circ (\iota\circ \zeta^{u,e})^{*} \circ (e_j^c)^* \Delta_{\Lm_k,\Lm_k} \\
&=& a_{\Sigma_{\bf{{\rm i}}}}^* \circ (\iota\circ \zeta^{u,e})^{*} c^{\delta_{j,k}}\Delta_{\Lm_k,\Lm_k}\\
&=& c^{\delta_{j,k}} a_{\Sigma_{\bf{{\rm i}}}}^* \Delta_{u\Lm_k,\Lm_k}=c^{\delta_{j,k}}A_{k_{\rm max}},
\end{eqnarray*}
where we use Proposition \ref{twistprop} in the third and fifth equality.
Thus, we obtain (\ref{claim5-2}).

Next, let us prove (\ref{claim5-3}). If 
$i\in \{1,2,\cdots,\tilde{r}\}\setminus \{i_1,i_2,\cdots,i_n\}$ then $u\Lm_i=\Lm_i$.
Since we know Proposition \ref{twistprop},
our claim (\ref{claim5-3}) follows from the following calculation: 
\begin{eqnarray*}
& &(a_{\Sigma_{\bf{{\rm i}}}}^{-1}\circ (\iota\circ \zeta^{u,e})^{-1} \circ e_j^c\circ (\iota\circ \zeta^{u,e}) \circ a_{\Sigma_{\bf{{\rm i}}}})^*A_{-i}\\
&=& a_{\Sigma_{\bf{{\rm i}}}}^* \circ (\iota\circ \zeta^{u,e})^{*} \circ (e_j^c)^*\circ
((\iota\circ \zeta^{u,e})^{-1})^*\circ(a_{\Sigma_{\bf{{\rm i}}}}^{-1})^*A_{-i}\\
&=& a_{\Sigma_{\bf{{\rm i}}}}^* \circ (\iota\circ \zeta^{u,e})^{*} \circ (e_j^c)^*\circ
((\iota\circ \zeta^{u,e})^{-1})^* \Delta_{\Lm_i,\Lm_i} \\
&=& a_{\Sigma_{\bf{{\rm i}}}}^* \circ (\iota\circ \zeta^{u,e})^{*} \circ (e_j^c)^* \circ ((\iota\circ \zeta^{u,e})^{-1})^*\Delta_{u\Lm_i,\Lm_i} \\
&=& a_{\Sigma_{\bf{{\rm i}}}}^* \circ (\iota\circ \zeta^{u,e})^{*} \circ (e_j^c)^* \Delta_{\Lm_i,\Lm_i} \\
&=&a_{\Sigma_{\bf{{\rm i}}}}^* \circ (\iota\circ \zeta^{u,e})^{*} \Delta_{\Lm_i,\Lm_i} \\
&=&a_{\Sigma_{\bf{{\rm i}}}}^* \Delta_{u\Lm_i,\Lm_i} = A_{-i}.
\end{eqnarray*}

We now turn to (\ref{claim5-1}) by induction on $n-k$.
First, let us consider the case $n-k=0$ so that $n=k$.
By Proposition \ref{compprop}(\ref{comp1}),
we have $(p\circ a_{\Sigma_{\bf{{\rm i}}}}^{-1}\circ (\iota\circ \zeta^{u,e})^{-1} \circ e_j^c \circ \iota\circ \zeta^{u,e} \circ a_{\Sigma_{\bf{{\rm i}}}}A)_{n}=\overline{X}_{n}$.
We set $F:=\{i\in I | i^+=n+1\}$.
Note that $\tilde{B}_{n,n^-}=-1$ and if $l\in I\setminus M\setminus \{n^-\}$ then $\tilde{B}_{n,l}=0$ by $l<l^+<n$ and Definition \ref{exmatdef} and the definition of the matrix $M$ in Theorem \ref{pgdef}.
Therefore, the definition of $p$ and (\ref{claim5-2}) imply
\begin{eqnarray*}
& &\overline{X}_{n}=(p\circ a_{\Sigma_{\bf{{\rm i}}}}^{-1}\circ (\iota\circ \zeta^{u,e})^{-1} \circ e_j^c \circ \iota\circ \zeta^{u,e} \circ a_{\Sigma_{\bf{{\rm i}}}}A)_{n}\\
&=& \prod_{l\in I}(a_{\Sigma_{\bf{{\rm i}}}}^{-1}\circ (\iota\circ \zeta^{u,e})^{-1} \circ e_j^c \circ \iota\circ \zeta^{u,e} \circ a_{\Sigma_{\bf{{\rm i}}}}A)_{l}^{\tilde{B}_{n,l}}\\
&=& \prod_{l\in I\setminus F} (a_{\Sigma_{\bf{{\rm i}}}}^{-1}\circ (\iota\circ \zeta^{u,e})^{-1} \circ e_j^c \circ \iota\circ \zeta^{u,e} \circ a_{\Sigma_{\bf{{\rm i}}}}A)_{l}^{\tilde{B}_{n,l}} \\
& & 
\prod_{l\in F} (a_{\Sigma_{\bf{{\rm i}}}}^{-1}\circ (\iota\circ \zeta^{u,e})^{-1} \circ e_j^c \circ \iota\circ \zeta^{u,e} \circ a_{\Sigma_{\bf{{\rm i}}}}A)_{l}^{\tilde{B}_{n,l}} \\
&=& (a_{\Sigma_{\bf{{\rm i}}}}^{-1}\circ (\iota\circ \zeta^{u,e})^{-1} \circ e_j^c \circ \iota\circ \zeta^{u,e} \circ a_{\Sigma_{\bf{{\rm i}}}}A)_{n^-}^{-1}c^{\tilde{B}_{n,j_{\rm max}}} \prod_{l\in F} A_l^{\tilde{B}_{n,l}}\\
&=& (a_{\Sigma_{\bf{{\rm i}}}}^{-1}\circ (\iota\circ \zeta^{u,e})^{-1} \circ e_j^c \circ \iota\circ \zeta^{u,e} \circ a_{\Sigma_{\bf{{\rm i}}}}A)_{n^-}^{-1}c^{\tilde{B}_{n,j_{\rm max}}} A_{n^-} p(A)_{n},
\end{eqnarray*}
which yields our claim $(a_{\Sigma_{\bf{{\rm i}}}}^{-1}\circ (\iota\circ \zeta^{u,e})^{-1} \circ e_j^c \circ \iota\circ \zeta^{u,e} \circ a_{\Sigma_{\bf{{\rm i}}}}A)_{n^-}=A_{n^-} c^{\tilde{B}_{n,j_{\rm max}}}\frac{p(A)_n}{\overline{X}_n}$ for $k=n$.

Next, let us consider the case $n-k>0$.
Using Proposition \ref{compprop} (\ref{comp1}),
we have $\overline{X}_{k}=(p\circ a_{\Sigma_{\bf{{\rm i}}}}^{-1}\circ  (\iota\circ \zeta^{u,e})^{-1} \circ e_j^c \circ (\iota\circ \zeta^{u,e}) \circ a_{\Sigma_{\bf{{\rm i}}}}A)_{k}$. Note that $\tilde{B}_{k,k^-}=-1$ and if $l\in \{-\tilde{r},\cdots,-1,1,\cdots,n\} \setminus F\setminus \{n^-,\cdots,(k+1)^-,k^-\}$ then $\tilde{B}_{k,l}=0$ by $l<l^+<k$ and Definition \ref{exmatdef} and the definition of the matrix $M$ in Theorem \ref{pgdef}. Hence the submatrix $(\tilde{B_{i,l}})_{k\leq i\leq n,\ l=(k+1)^-,\cdots,n^-,j_{\rm max}}$ of $\tilde{B}_{\textbf{i}}$ is as follows:
\[
  \left(
    \begin{array}{ccccc}
      \tilde{B}_{k,(k+1)^-} & \tilde{B}_{k,(k+2)^-} & \ldots & \tilde{B}_{k,n^-} & \tilde{B}_{k,j_{\rm max}} \\
      -1 & \tilde{B}_{k+1,(k+2)^-} & \ldots & \tilde{B}_{k+1,n^-} & \tilde{B}_{k+1,j_{\rm max}}\\
      0 & -1 & \ldots & \tilde{B}_{k+2,n^-} & \tilde{B}_{k+2,j_{\rm max}}\\      
      \vdots & \vdots & \ddots & \vdots & \vdots \\
      0 & 0 & \ldots & -1 & \tilde{B}_{n,j_{\rm max}}
    \end{array}
  \right)
\]

In conjunction with the definition of $p$,
(\ref{claim5-2}) and induction hypothesis, we obtain
\begin{eqnarray*}
& &\overline{X}_{k}=(p\circ a_{\Sigma_{\bf{{\rm i}}}}^{-1}\circ (\iota\circ \zeta^{u,e})^{-1} \circ e_j^c \circ (\iota\circ \zeta^{u,e}) \circ a_{\Sigma_{\bf{{\rm i}}}}A)_{k}\\
&=& \prod_{l\in \{-\tilde{r},\cdots,-1,1,\cdots,n\}} (a_{\Sigma_{\bf{{\rm i}}}}^{-1}\circ (\iota\circ \zeta^{u,e})^{-1} \circ e_j^c \circ (\iota\circ \zeta^{u,e}) \circ a_{\Sigma_{\bf{{\rm i}}}}A)_l^{\tilde{B}_{k,l}}\\
&=&
(a_{\Sigma_{\bf{{\rm i}}}}^{-1}\circ (\iota\circ \zeta^{u,e})^{-1} \circ e_j^c \circ (\iota\circ \zeta^{u,e}) \circ a_{\Sigma_{\bf{{\rm i}}}}A)_{k^-}^{-1}\\
& &\cdot 
\prod_{l\in \{-\tilde{r},\cdots,-1,1,\cdots,n\}\setminus F\setminus \{n^-,\cdots,(k+1)^-,k^-\}}
(a_{\Sigma_{\bf{{\rm i}}}}^{-1}\circ (\iota\circ \zeta^{u,e})^{-1} \circ e_j^c \circ (\iota\circ \zeta^{u,e}) \circ a_{\Sigma_{\bf{{\rm i}}}}A)_l^{\tilde{B}_{k,l}}\\
& &\cdot
\prod_{l\in \{n^-,\cdots,(k+1)^-\}}
(a_{\Sigma_{\bf{{\rm i}}}}^{-1}\circ (\iota\circ \zeta^{u,e})^{-1} \circ e_j^c \circ (\iota\circ \zeta^{u,e}) \circ a_{\Sigma_{\bf{{\rm i}}}}A)_l^{\tilde{B}_{k,l}}\\
& &\cdot \prod_{l\in F}
(a_{\Sigma_{\bf{{\rm i}}}}^{-1}\circ (\iota\circ \zeta^{u,e})^{-1} \circ e_j^c \circ (\iota\circ \zeta^{u,e}) \circ a_{\Sigma_{\bf{{\rm i}}}}A)_l^{\tilde{B}_{k,l}}\\
&=& (a_{\Sigma_{\bf{{\rm i}}}}^{-1}\circ (\iota\circ \zeta^{u,e})^{-1} \circ e_j^c \circ (\iota\circ \zeta^{u,e}) \circ a_{\Sigma_{\bf{{\rm i}}}}A)_{k^-}^{-1}\\
& &\cdot
\prod_{l\in \{n^-,\cdots,(k+1)^-\}}
(a_{\Sigma_{\bf{{\rm i}}}}^{-1}\circ (\iota\circ \zeta^{u,e})^{-1} \circ e_j^c \circ (\iota\circ \zeta^{u,e}) \circ a_{\Sigma_{\bf{{\rm i}}}}A)_l^{\tilde{B}_{k,l}}\\
& &\cdot c^{\tilde{B}_{k,j_{\rm max}}} \prod_{l\in F}
A_l^{\tilde{B}_{k,l}}\\
&=& (a_{\Sigma_{\bf{{\rm i}}}}^{-1}\circ (\iota\circ \zeta^{u,e})^{-1} \circ e_j^c \circ (\iota\circ \zeta^{u,e}) \circ a_{\Sigma_{\bf{{\rm i}}}}A)_{k^-}^{-1}\cdot c^{\tilde{B}_{k,j_{\rm max}}} \prod_{l\in F}A_l^{\tilde{B}_{k,l}}\\
& &\cdot
\prod_{l\in \{n,\cdots,k+1\}}
\left(A_{l^-} c^{D_{\{l,l+1,\cdots,n\},\{(l+1)^-,(l+2)^-,\cdots,n^-,j_{\rm max}\}}}
\left(\frac{p(A)_{l}}{\overline{X}_{l}}\right)
\left(\frac{p(A)_{l+1}}{\overline{X}_{l+1}}\right)^{D_{l,(l+1)^-}} \right. \\
& &\left.\cdot 
\left( \frac{p(A)_{l+2}}{\overline{X}_{l+2}}\right)^{D_{\{l,l+1\},\{(l+1)^-,(l+2)^-\}}}
\cdots
\left(\frac{p(A)_{n}}{\overline{X}_{n}}\right)^{D_{\{l,l+1,\cdots,n-1\},\{(l+1)^-,(l+2)^-,\cdots,n^-\}}} \right)^{\tilde{B}_{k,l^-}}\\
&=& (a_{\Sigma_{\bf{{\rm i}}}}^{-1}\circ (\iota\circ \zeta^{u,e})^{-1} \circ e_j^c \circ (\iota\circ \zeta^{u,e}) \circ a_{\Sigma_{\bf{{\rm i}}}}A)_{k^-}^{-1}\cdot A_{k^-} p(A)_{k}\\
& &\cdot
c^{\tilde{B}_{k,j_{\rm max}}+\sum_{l=n,n-1,\cdots,k+1}D_{\{l,l+1,\cdots,n\},\{(l+1)^-,\cdots,n^-,j_{\rm max}\}}\tilde{B}_{k,l^-}}\\
& &\cdot
\prod_{l\in \{n,\cdots,k+1\}}
\left(\frac{p(A)_{l}}{\overline{X}_{l}}\right)^{\tilde{B}_{k,l^-}+\sum_{s=l-1,l-2,\cdots,k+1} \tilde{B}_{k,s^-}D_{\{s,\cdots,l-2,l-1\},\{(s+1)^-,\cdots,(l-1)^-,l^-\}}} \\
&=& (a_{\Sigma_{\bf{{\rm i}}}}^{-1}\circ (\iota\circ \zeta^{u,e})^{-1} \circ e_j^c \circ (\iota\circ \zeta^{u,e}) \circ a_{\Sigma_{\bf{{\rm i}}}}A)_{k^-}^{-1}\cdot A_{k^-} p(A)_{k}\\
& &\cdot
c^{D_{\{k,k+1,\cdots,n\},\{(k+1)^-,\cdots,n^-,j_{\rm max}\}}}\\
& &\cdot
\prod_{l\in \{n,\cdots,k+1\}}
\left(\frac{p(A)_{l}}{\overline{X}_{l}}\right)^{D_{\{k,k+1,\cdots,l-1\},\{(k+1)^-,\cdots,l^-\}}},
\end{eqnarray*}
which yields our claim (\ref{claim5-1}). \qed

\begin{thm}\label{thm5}

Let ${\rm \bf{i}}=(i_1,\cdots,i_n)$ be a reduced word of $u\in W$,
$A=(A_{-\tilde{r}},\cdots,A_{-1},A_1,\cdots,A_n)\in \mathcal{A}_{\Sigma_{\rm\bf{i}}}$, $M$ be the matrix in Theorem \ref{pgdef}, and $p=p_M:\mathcal{A}_{\Sigma_{\rm\bf{i}}}\rightarrow \mathcal{X}_{\Sigma_{\rm\bf{i}}}$ be the map in Proposition \ref{map-p}. For $j\in[1,r]$, we set $\{K_1,K_2,\cdots,K_l\}:=\{K|1\leq K\leq n ,\ i_K=j\}$ $(K_1<\cdots<K_l)$.
\[
\gamma_j\circ (\iota\circ \zeta^{u,e})\circ a_{\Sigma_{\rm\bf{i}}}(A)=
p(A)_{-j}p(A)_{K_1}\cdots p(A)_{K_l},
\]
\[
\varepsilon_j\circ (\iota\circ \zeta^{u,e})\circ a_{\Sigma_{\rm\bf{i}}}(A)=
\left(\sum^{l-1}_{s=0}p(A)_{K_{s+1}}p(A)_{K_{s+2}}\cdots p(A)_{K_l}\right)^{-1}.
\]
\end{thm}

\nd
{\it Proof.}

It follows from Proposition \ref{compprop} (\ref{comp2}) that
$\varepsilon_j\circ x_{\Sigma_{\rm\bf{i}}}\circ p=\varepsilon_j\circ (\iota\circ \zeta^{u,e})\circ a_{\Sigma_{\rm\bf{i}}}$ and $\gamma_j\circ x_{\Sigma_{\rm\bf{i}}}\circ p=\gamma_j\circ (\iota\circ \zeta^{u,e})\circ a_{\Sigma_{\rm\bf{i}}}$.
Our claims follow from Theorem \ref{thm1}. \qed

\vspace{3mm}

Theorem \ref{thm4} and Theorem \ref{thm5} imply that the map 
$(\iota\circ \zeta^{u,e})\circ a_{\Sigma_{\rm\bf{i}}}:\mathcal{A}_{\Sigma_{\rm \bf{i}}}\rightarrow G^{u,e}$
is a positive structure on the geometric crystal
$(G^{u,e}, \{e_i\}_{i\in [1,r]},\{\gamma_i\}_{i\in [1,r]},\{\varepsilon_i\}_{i\in [1,r]})$. 
Let $\Sigma$ be a seed obtained from $\Sigma_{\rm\bf{i}}$ by an iteration
of mutations. Then the corresponding birational map 
$\mathcal{A}_{\Sigma_{\rm \bf{i}}}\rightarrow \mathcal{A}_{\Sigma}$
is positive. Hence, we obtain the following corollary.:

\begin{cor}\label{corA}
The map $(\iota\circ \zeta^{u,e})\circ a_{\Sigma}:\mathcal{A}_{\Sigma}\rightarrow G^{u,e}$ is a
positive structure on the geometric crystal
$(G^{u,e}, \{e_i\}_{i\in [1,r]},\{\gamma_i\}_{i\in [1,r]},\{\varepsilon_i\}_{i\in [1,r]})$. 
Applying Theorem \ref{UDthm}, we obtain a crystal
$(X_{*}(\mathcal{A}_{\Sigma}), \{\tilde{e}_i\}_{i\in [1,r]},\{\tilde{\gamma}_i\}_{i\in [1,r]},\{\tilde{\varepsilon}_i\}_{i\in [1,r]})$.
It induces a glued crystal structure on 
the set of $\mathbb{Z}$-valued points $\mathcal{A}_{|\Sigma_{\rm \bf{i}}|}(\mathbb{Z}^T)$.
\end{cor}

The glued crystals in Corollary \ref{corX}, \ref{corA} seem to be related to
Fock-Goncharov conjecture
as seen in \ref{zvalue}.

\section{Type A-case}\label{7sec}

In the rest of article, we set $G=SL_{r+1}(\mathbb{C})$ and consider the cell $G^{w_0,e}$ with the longest element $w_0\in W$. 
Let $n:=\frac{r(r+1)}{2}$, $I=\{-r,\cdots,-1,1,2,\cdots,n\}$ and $\textbf{i}_0$ be the following reduced word of $w_0$:
\[
\textbf{i}_0=(\underbrace{1,2,\cdots,r-2,r-1,r}_{1\ {\rm st\ cycle}},\underbrace{1,2,\cdots,r-2,r-1}_{2\ {\rm nd\ cycle}},\underbrace{1,2,\cdots,r-2}_{3\ {\rm rd\ cycle}},\cdots,\underbrace{1,2}_{r-1\ {\rm th\ cycle}},1)
\]
and $i_k$ be the $k$-th index of $\textbf{i}_0$ from the left.

In this section, we shall calculate explicit forms of
$a_{\Sigma_{\bf{{\rm i}}_0}}^{-1}\circ (\iota\circ\zeta^{w_0,e})^{-1} \circ e_j^c \circ  (\iota\circ\zeta^{w_0,e}) \circ a_{\Sigma_{\bf{{\rm i}}_0}}$
by a direct calculation. 
In the last subsection,
we will verify the explicit forms coincide with those of Theorem \ref{thm4}. 

\subsection{Fundamental representation of Type ${\rm A}_r$}

First, we review the fundamental representations of the complex simple Lie algebras $\ge$ of type ${\rm A}_r$ \cite{KN:1994, N}. 
Let $\frg=\frs\frl(r+1,\mathbb{C})$ be the simple Lie algebra of type ${\rm A}_r$. 
The Cartan matrix $A=(a_{i,j})_{i,j\in \{1,2,\cdots,r\}}$ of $\frg$ is as follows:
\[a_{i,j}=
\begin{cases}
2 & {\rm if}\ i=j, \\
-1 & {\rm if}\ |i-j|=1, \\
0 & {\rm otherwise.}  
\end{cases}
\]
For $\frg=\lan \frh,e_i,f_i(i\in \{1,2,\cdots,r\})\ran$, 
let us describe the vector representation 
$V(\Lm_1)$. Set ${\mathbf B}^{(r)}:=
\{v_i|\ i=1,2,\cd,r+1\}$ and define 
$V(\Lm_1):=\bigoplus_{v\in{\mathbf B}^{(r)}}\bbC v$. The weights of $v_i$ $(i=1,\cd,r+1)$ are given by ${\rm wt}(v_i)=\Lm_i-\Lm_{i-1}$, where $\Lm_0=\Lm_{r+1}=0$. We define the $\frg$-action on $V(\Lm_1)$ as follows:
\begin{eqnarray}
&& h v_j=\lan h,{\rm wt}(v_j)\ran v_j\ \ (h\in P^*,\ j\in \{1,2,\cdots,r+1\}), \\
&&f_iv_i=v_{i+1},\q
e_iv_{i+1}=v_i \q(1\leq i\leq r),\label{A-f1}
\end{eqnarray}
and the other actions are trivial.

Let $\Lm_i$ be the $i$-th fundamental weight of type ${\rm A}_r$.
As is well-known that the fundamental representation 
$V(\Lm_i)$ $(1\leq i\leq r)$
is embedded in $\wedge^i V(\Lm_1)$
with multiplicity free.
The explicit form of the highest (resp. lowest) weight 
vector $v_{\Lm_i}$ (resp. $u_{\Lm_i}$)
of $V(\Lm_i)$ is realized in 
$\wedge^i V(\Lm_1)$ as follows:
\begin{equation}\label{A-h-l}
v_{\Lm_i}=v_1\wedge v_2\wedge\cdots\wedge v_i, \qq
u_{\Lm_i}=v_{r-i+2}\wedge v_{r-i+3}\wedge\cdots\wedge v_{r+1}.
\end{equation}

It is known that if $1\leq j_1<\cdots<j_s\leq r+1$, $1\leq l_1<\cdots<l_s\leq r+1$ and $x\in G=SL_{r+1}(\mathbb{C})$ then
the value of bilinear form in \ref{bilin} on $x$ 
\[
\lan v_{j_1}\wedge v_{j_2} \wedge \cdots \wedge v_{j_s}, x v_{l_1}\wedge v_{l_2} \wedge \cdots \wedge v_{l_s}  \ran
\]
is equal to an ordinary minor $D_{\{j_1,\cdots,j_s\},\{l_1,\cdots,l_s\}}(x)$.

\subsection{Geometric crystal action on cluster $\mathcal{A}$-varieties for type ${\rm A}_r$ case}


\begin{lem}\label{minordesc}
The rational maps $\gamma_j$, $\varepsilon_j$ : $G^{w_0,e}\rightarrow \mathbb{C}^{\times}$ can be described as
\[
\varepsilon_j=\frac{D_{j+1,j+1}}{D_{j+1,j}},\ \ 
\gamma_j=\frac{D_{j,j}}{D_{j+1,j+1}}.
\]
\end{lem}

\nd
{\it Proof.}

Our claim follows by (\ref{minorexp}) and (\ref{A-h-l}). \qed

\begin{lem}\label{GC-lem}
We suppose that $1\leq j_1<\cdots<j_s\leq r+1$ and $j\in\{1,2,\cdots,r\}$. 
If there exists $i\in\{1,2,\cdots,s\}$ such that
$j=j_i$ and $j_{i+1}>j_i+1$ (we set $j_{s+1}=r+2$) then
\[
(e_j^c)^* D_{\{j_1,\cdots,j_s\},\{1,\cdots,s\}}
=D_{\{j_1,\cdots,j_s\},\{1,\cdots,s\}}+(c-1) \frac{D_{\{1,\cdots,j\},\{1,\cdots,j\}}}{D_{\{1,\cdots,j-1,j+1\},\{1,\cdots,j\}}}D_{\{j_1,\cdots,j_{i-1},j_i+1,j_{i+1},\cdots, j_s\},\{1,\cdots,s\}},
\]
otherwise,
\[
(e_j^c)^* D_{\{j_1,\cdots,j_s\},\{1,\cdots,s\}}
=D_{\{j_1,\cdots,j_s\},\{1,\cdots,s\}}.
\]
\end{lem}

\nd
{\it Proof.}

Because of Lemma \ref{minordesc}, we obtain
\[
\varphi_j=\varepsilon_j\gamma_j =\frac{D_{j,j}}{D_{j+1,j}} =
 \frac{D_{\{1,\cdots,j\},\{1,\cdots,j\}}}{D_{\{1,\cdots,j-1,j+1\},\{1,\cdots,j\}}}.
\]

We can calculate
\begin{eqnarray*}
& &(e_j^c)^* D_{\{j_1,\cdots,j_s\},\{1,\cdots,s\}}(x)\\
&=& \lan v_{j_1}\wedge v_{j_2} \wedge \cdots \wedge v_{j_s},
e_j^c(x) v_{1}\wedge v_{2} \wedge \cdots \wedge v_{s}  \ran \\
&=& \lan v_{j_1}\wedge v_{j_2} \wedge \cdots \wedge v_{j_s},
x_j((c-1)\varphi_j(x))x x_j((c^{-1}-1)\varepsilon_i(x)) v_{1}\wedge v_{2} \wedge \cdots \wedge v_{s}  \ran \\
&=& \lan y_j((c-1)\varphi_j(x)) v_{j_1}\wedge v_{j_2} \wedge \cdots \wedge v_{j_s},
x v_{1}\wedge v_{2} \wedge \cdots \wedge v_{s}  \ran. \\
\end{eqnarray*}

\vspace{-5mm}

If there exists $i\in\{1,2,\cdots,s\}$ such that $j=j_i$ and $j_{i+1}>j_i+1$ then 
\begin{eqnarray*}
& &
(e_j^c)^* D_{\{j_1,\cdots,j_s\},\{1,\cdots,s\}}(x)\\
&=&
\lan (1+(c-1)\varphi_j(x)f_j ) v_{j_1}\wedge v_{j_2} \wedge \cdots \wedge v_{j_s},
x v_{1}\wedge v_{2} \wedge \cdots \wedge v_{s}  \ran \\
&=&
\lan v_{j_1}\wedge v_{j_2} \wedge \cdots \wedge v_{j_s}, x v_{1}\wedge v_{2} \wedge \cdots \wedge v_{s}  \ran \\
& &+(c-1)\varphi_j(x) \lan v_{j_1}\wedge \cdots \wedge v_{j_{i-1}}\wedge v_{j_i+1} \wedge \cdots \wedge v_{j_s},
x v_{1}\wedge v_{2} \wedge \cdots \wedge v_{s}  \ran \\
&=& D_{\{j_1,\cdots,j_s\},\{1,\cdots,s\}}(x) + (c-1)\varphi_j(x) D_{\{j_1,\cdots,j_{i-1},j_i+1,j_{i+1},\cdots,j_s\},\{1,\cdots,s\}}(x) \\
&=& D_{\{j_1,\cdots,j_s\},\{1,\cdots,s\}}(x) + (c-1)\frac{D_{\{1,\cdots,j\},\{1,\cdots,j\}}}{D_{\{1,\cdots,j-1,j+1\},\{1,\cdots,j\}}} D_{\{j_1,\cdots,j_{i-1},j_i+1,j_{i+1},\cdots,j_s\},\{1,\cdots,s\}}(x).
\end{eqnarray*}
Otherwise, we have
$(e_j^c)^* D_{\{j_1,\cdots,j_s\},\{1,\cdots,s\}}(x)
= \lan v_{j_1}\wedge v_{j_2} \wedge \cdots \wedge v_{j_s}, x v_{1}\wedge v_{2} \wedge \cdots \wedge v_{s}  \ran
= D_{\{j_1,\cdots,j_s\},\{1,\cdots,s\}}(x)$. \qed

\vspace{2mm}

Recall that a regular map
$p_M : \mathcal{A}_{\Sigma_{\textbf{i}_0}}\rightarrow \mathcal{X}_{\Sigma_{\textbf{i}_0}}=\mathcal{X}_{\textbf{i}_0}$
is defined as
\[
p^* (X_i)=\prod_{k\in I}A^{\tilde{B}_{i,k}}_k,
\]
where $(X_i)_{i\in I}$ is a coordinate of $\mathcal{X}_{\Sigma_{\textbf{i}_0}}$.
In what follows, we denote $p_M$ by $p$.


We write a coordinate $A=(A_{-r},\cdots,A_{-1},A_1,\cdots,A_n)\in \mathcal{A}_{\Sigma_{\textbf{i}_0}}$
as
\[
(A_{-r},\cdots,A_{-1},A_{1,1},A_{1,2},\cdots,A_{1,r},A_{2,1},A_{2,2},\cdots,A_{2,r-1},\cdots,A_{r-1,1},A_{r-1,2},A_{r,1}).
\]
Note that $(a_{\Sigma_{\bf{{\rm i}}_0}}^*)^{-1} A_{m,d}=\lan v_{m+1}\wedge v_{m+2}\wedge \cdots \wedge v_{m+d}, \cdot v_1\wedge v_2\wedge \cdots \wedge v_d \ran$ for $m\in\{1,2,\cdots,r\}$ and $d\in\{1,2,\cdots,r-m+1\}$.

We also set
\[
(P_{-r},\cdots,P_{-1},P_{1,1},P_{1,2},\cdots,P_{1,r},P_{2,1},P_{2,2},\cdots,P_{2,r-1},\cdots,P_{r-1,1},P_{r-1,2},P_{r,1})=p(A).
\]
\begin{thm}\label{thm2}
For $m\in\{1,2,\cdots,r\}$ and $d\in\{1,2,\cdots,r-m+1\}$, we obtain

\vspace{3mm}

\hspace{-5mm}
$(a_{\Sigma_{\bf{{\rm i}}_0}}^{-1}\circ (\iota\circ\zeta^{w_0,e})^{-1} \circ e_j^c \circ  (\iota\circ\zeta^{w_0,e}) \circ a_{\Sigma_{\bf{{\rm i}}_0}}A)_{m,d}$
\[=
\begin{cases}
 A_{m,d} \frac{c(P_{1,d}P_{2,d}\cdots P_{r-d,d}+P_{2,d}\cdots P_{r-d,d}+\cdots+ P_{m,d}\cdots P_{r-d,d})+P_{m+1,d}\cdots P_{r-d,d}+\cdots +P_{r-d,d}+1}{P_{1,d}P_{2,d}\cdots P_{r-d,d}+P_{2,d}\cdots P_{r-d,d}+\cdots +P_{r-d,d}+1} & {\rm if}\ j=d,\\
A_{m,d} & {\rm if}\ j\neq d.
\end{cases}
\]
\end{thm}

\nd
{\it Proof.}

By Proposition \ref{twistprop} and Lemma \ref{GC-lem}, we get
\begin{eqnarray*}
& &A_{m,d} \overset{(a_{\Sigma_{\bf{{\rm i}}_0}}^{-1})^*}{\longrightarrow}
\lan v_{m+1}\wedge v_{m+2}\wedge \cdots \wedge v_{m+d}, \cdot v_1\wedge v_2\wedge \cdots \wedge v_d \ran \\
& &
=\frac{\lan v_{m+1}\wedge v_{m+2}\wedge \cdots \wedge v_{r+1}, \cdot v_1\wedge v_2\wedge \cdots \wedge v_d 
\wedge v_{m+d+1}\wedge \cdots \wedge v_{r+1}\ran}{\lan v_{m+d+1}\wedge v_{m+d+2}\wedge \cdots \wedge v_{r+1}, \cdot
v_{m+d+1}\wedge \cdots \wedge v_{r+1}\ran} \\
& &
=\frac{\lan \ovl{w_0}v_{1}\wedge v_{2}\wedge \cdots \wedge v_{r-m+1}, \cdot v_1\wedge v_2\wedge \cdots \wedge v_d 
\wedge v_{m+d+1}\wedge \cdots \wedge v_{r+1}\ran}{\lan \ovl{w_0}v_{1}\wedge v_{2}\wedge \cdots \wedge v_{r-m-d+1}, \cdot
v_{m+d+1}\wedge \cdots \wedge v_{r+1}\ran}
\\
& &
\overset{((\iota\circ \zeta^{w_0,e})^{-1})^*}{\longrightarrow}
\frac{\lan v_1\wedge v_2\wedge \cdots \wedge v_d 
\wedge v_{m+d+1}\wedge \cdots \wedge v_{r+1}, \cdot v_{1}\wedge v_{2}\wedge \cdots \wedge v_{r-m+1}\ran}{\lan v_{m+d+1}\wedge \cdots \wedge v_{r+1}, \cdot v_{1}\wedge v_{2}\wedge \cdots \wedge v_{r-m-d+1}\ran} \\
& &
=\frac{D_{\{1,\cdots,d,m+d+1,\cdots,r+1\},\{1,2,\cdots,r-m+1\}}}{D_{\{m+d+1,\cdots,r+1\},\{1,2,\cdots,r-m-d+1\}}}\\
& & \overset{(e_d^c)^*}{\longrightarrow}
\frac{1}{D_{\{m+d+1,\cdots,r+1\},\{1,2,\cdots,r-m-d+1\}}} 
\left(D_{\{1,\cdots,d,m+d+1,\cdots,r+1\},\{1,2,\cdots,r-m+1\}} \right. \\
& &
\left.+(c-1)\frac{D_{\{1,\cdots,d\},\{1,\cdots,d\}}}{D_{\{1,\cdots,d-1,d+1\},\{1,\cdots,d\}}} D_{\{1,\cdots,d-1,d+1,m+d-1,\cdots,r+1\}, \{1,2,\cdots,r-m+1\}} \right)\\
& & \overset{(\iota\circ \zeta^{w_0,e})^*}{\longrightarrow}
\frac{1}{D_{\{m+d+1,\cdots,r+1\},\{m+d+1,\cdots,r+1\}}} 
\left(D_{\{m+1,\cdots,r+1\},\{1,2,\cdots,d,m+d+1,\cdots,r+1\}} \right. \\
& &
\left.+(c-1)\frac{D_{\{r-d+2,\cdots,d\},\{1,\cdots,d\}}}{D_{\{r-d+2,\cdots,r+1\},\{1,\cdots,d-1,d+1\}}} D_{\{m+1,\cdots,r+1\}, \{1,2,\cdots,d-1,d+1,m+d+1,\cdots,r-m+1\}} \right).
\end{eqnarray*}

Hence, we need to show that
\begin{eqnarray}
& &(a_{\Sigma_{\bf{{\rm i}}_0}})^* \frac{1}{D_{\{m+d+1,\cdots,r+1\},\{m+d+1,\cdots,r+1\}}} 
\left(D_{\{m+1,\cdots,r+1\},\{1,2,\cdots,d,m+d+1,\cdots,r+1\}} \right. \nonumber\\
& &
\left.+(c-1)\frac{D_{\{r-d+2,\cdots,r+1\},\{1,\cdots,d\}}}{D_{\{r-d+2,\cdots,r+1\},\{1,\cdots,d-1,d+1\}}} D_{\{m+1,\cdots,r+1\}, \{1,2,\cdots,d-1,d+1,m+d+1,\cdots,r-m+1\}} \right) \label{pr-1} \\
&=& A_{m,d} \frac{c(P_{1,d}P_{2,d}\cdots P_{r-d,d}+P_{2,d}\cdots P_{r-d,d}+\cdots+ P_{m,d}\cdots P_{r-d,d})+P_{m+1,d}\cdots P_{r-d,d}+\cdots +P_{r-d,d}+1}{P_{1,d}P_{2,d}\cdots P_{r-d,d}+P_{2,d}\cdots P_{r-d,d}+\cdots +P_{r-d,d}+1}. \nonumber
\end{eqnarray}

Let us take an element of the open subset $\ovl{\mathbb{B}}^-_{w_0}$ in (\ref{anopen})
\[
x=a y_{r}(t_{1,r}) \cdots y_{1}(t_{1,1})y_{r}(t_{2,r}) \cdots y_{2}(t_{2,2})\cdots y_{r}(t_{r-1,r})y_{r-1}(t_{r-1,r-1}) y_{r}(t_{r,r})
\in \ovl{\mathbb{B}}^-_{w_0}\subset G^{w_0,e}
\]
with $a={\rm diag}(a_1,\cdots,a_r,a_{r+1})\in H$ and $t_{s,i}\in\mathbb{C}^{\times}$. We get
\begin{eqnarray}
& & \frac{1}{D_{\{m+d+1,\cdots,r+1\},\{m+d+1,\cdots,r+1\}}} 
\left(D_{\{m+1,\cdots,r+1\},\{1,2,\cdots,d,m+d+1,\cdots,r+1\}} \right. \nonumber\\
& &
\left.+(c-1)\frac{D_{\{r-d+2,\cdots,r+1\},\{1,\cdots,d\}}}{D_{\{r-d+2,\cdots,r+1\},\{1,\cdots,d-1,d+1\}}} D_{\{m+1,\cdots,r+1\}, \{1,2,\cdots,d-1,d+1,m+d+1,\cdots,r-m+1\}} \right) (x) \nonumber \\
&=&a_{m+1}\cdots a_{m+d} \left(\prod^{d}_{i=1} t_{i,i}t_{i,i+1}\cdots t_{i,m+i-1} 
 +(c-1)\prod^{d}_{i=1} t_{i,i}t_{i,i+1}\cdots t_{i,r-d+i} \right. \nonumber\\
& & \left(\prod^{d-1}_{i=1} t_{i,i}t_{i,i+1}\cdots t_{i,r-d+i}(\sum^{r}_{i=d}t_{d+1,d+1}t_{d+1,d+2}\cdots t_{d+1,i}t_{d,i+1}\cdots t_{d,r})\right)^{-1} \nonumber \\
& & \left. \left(\prod^{d-1}_{i=1} t_{i,i}t_{i,i+1}\cdots t_{i,m+i-1}(\sum^{m+d-1}_{i=d}t_{d+1,d+1}t_{d+1,d+2}\cdots t_{d+1,i}t_{d,i+1}\cdots t_{d,m+d-1})\right) \right) \nonumber\\
&=&\ a_{m+1}\cdots a_{m+d}\prod^{d}_{i=1} t_{i,i}t_{i,i+1}\cdots t_{i,m+i-1} \left(1
 +(c-1) t_{d,m+d}t_{d,m+d+1}\cdots t_{d,r} \right. \nonumber \\
& &\left. \frac{\sum^{m+d-1}_{i=d}t_{d+1,d+1}t_{d+1,d+2}\cdots t_{d+1,i}t_{d,i+1}\cdots t_{d,m+d-1}}{\sum^{r}_{i=d}t_{d+1,d+1}t_{d+1,d+2}\cdots t_{d+1,i}t_{d,i+1}\cdots t_{d,r}} \right) \nonumber\\
&=&\ a_{m+1}\cdots a_{m+d}\prod^{d}_{i=1} t_{i,i}t_{i,i+1}\cdots t_{i,m+i-1} \label{pr-2} \\
& &\frac{\sum^{r}_{i=m+d} t_{d+1,d+1}t_{d+1,d+2}\cdots t_{d+1,i}t_{d,i+1}\cdots t_{d,r} + c\sum^{m+d-1}_{i=d} t_{d+1,d+1}\cdots t_{d+1,i}t_{d,i+1}\cdots t_{d,r}}{\sum^{r}_{i=d} t_{d+1,d+1}t_{d+1,d+2}\cdots t_{d+1,i}t_{d,i+1}\cdots t_{d,r}}. \nonumber
\end{eqnarray}
Next, let us calculate 
\[
(a_{\Sigma_{\bf{{\rm i}}_0}}^{-1})^*
A_{m,d} \frac{c(P_{1,d}\cdots P_{r-d,d}+\cdots+ P_{m,d}\cdots P_{r-d,d})+P_{m+1,d}\cdots P_{r-d,d}+\cdots +P_{r-d,d}+1}{P_{1,d}P_{2,d}\cdots P_{r-d,d}+P_{2,d}\cdots P_{r-d,d}+\cdots +P_{r-d,d}+1}(x).
\]
The construction of $\tilde{B}_{\rm \bf{i}_0}$ means
$P_{s,d}=\frac{A_{s+1,d}A_{s,d-1}A_{s-1,d+1}}{A_{s+1,d-1}A_{s,d+1}A_{s-1,d}}$ ($s=1,2,\cdots,r-d$).
We obtain
\[
(a_{\Sigma_{\bf{{\rm i}}_0}}^{-1})^* A_{s+1,d}(x)=D_{\{s+2,\cdots,s+d+1\},\{1,\cdots,d\}}(x)
=a_{s+2}\cdots a_{s+d+1} \prod^{d}_{i=1} t_{i,i}t_{i,i+1}\cdots t_{i,s+i},
\]
\[
(a_{\Sigma_{\bf{{\rm i}}_0}}^{-1})^* A_{s,d-1}(x)=D_{\{s+1,\cdots,s+d-1\},\{1,\cdots,d-1\}}(x)
=a_{s+1}\cdots a_{s+d-1} \prod^{d-1}_{i=1} t_{i,i}t_{i,i+1}\cdots t_{i,s+i-1},
\]
\[
(a_{\Sigma_{\bf{{\rm i}}_0}}^{-1})^* A_{s-1,d+1}(x)=D_{\{s,\cdots,s+d\},\{1,\cdots,d+1\}}(x)
=a_{s}\cdots a_{s+d} \prod^{d+1}_{i=1} t_{i,i}t_{i,i+1}\cdots t_{i,s+i-2},
\]
\[
(a_{\Sigma_{\bf{{\rm i}}_0}}^{-1})^* A_{s+1,d-1}(x)=D_{\{s+2,\cdots,s+d\},\{1,\cdots,d-1\}}(x)
=a_{s+2}\cdots a_{s+d} \prod^{d-1}_{i=1} t_{i,i}t_{i,i+1}\cdots t_{i,s+i},
\]
\[
(a_{\Sigma_{\bf{{\rm i}}_0}}^{-1})^* A_{s,d+1}(x)=D_{\{s+1,\cdots,s+d+1\},\{1,\cdots,d+1\}}(x)
=a_{s+1}\cdots a_{s+d+1} \prod^{d+1}_{i=1} t_{i,i}t_{i,i+1}\cdots t_{i,s+i-1},
\]
\[
(a_{\Sigma_{\bf{{\rm i}}_0}}^{-1})^* A_{s-1,d}(x)=D_{\{s,\cdots,s+d-1\},\{1,\cdots,d\}}(x)
=a_{s}\cdots a_{s+d-1} \prod^{d}_{i=1} t_{i,i}t_{i,i+1}\cdots t_{i,s+i-2}.
\]
Consequently, we have $(a_{\Sigma_{\bf{{\rm i}}_0}}^{-1})^*P_{s,d}(x)=\frac{t_{d,s+d}}{t_{d+1,s+d}}$. In conjunction with
$(a_{\Sigma_{\bf{{\rm i}}_0}}^{-1})^* A_{m,d}(x)=D_{\{m+1,\cdots,m+d\},\{1,\cdots,d\}}(x)
=a_{m+1}\cdots a_{m+d}\prod^d_{i=1}t_{i,i}t_{i,i+1}\cdots t_{i,s+i-1}$, one obtain
\begin{eqnarray}
& &(a_{\Sigma_{\bf{{\rm i}}_0}}^{-1})^*
A_{m,d} \frac{c(P_{1,d}\cdots P_{r-d,d}+\cdots+ P_{m,d}\cdots P_{r-d,d})+P_{m+1,d}\cdots P_{r-d,d}+\cdots +P_{r-d,d}+1}{P_{1,d}P_{2,d}\cdots P_{r-d,d}+P_{2,d}\cdots P_{r-d,d}+\cdots +P_{r-d,d}+1}(x) \nonumber \\
&=&a_{m+1}\cdots a_{m+d}\prod^d_{i=1}t_{i,i}t_{i,i+1}\cdots t_{i,s+i-1}\nonumber \\
& & \frac{c (\sum^{m}_{s=1} \frac{t_{d,s+d}t_{d,s+d+1}\cdots t_{d,r}}{t_{d+1,s+d}t_{d+1,s+d+1}\cdots t_{d+1,r}})
 + \sum^{r-d+1}_{s=m+1} \frac{t_{d,s+d}t_{d,s+d+1}\cdots t_{d,r}}{t_{d+1,s+d}t_{d+1,s+d+1}\cdots t_{d+1,r}}}
{\sum^{r-d+1}_{s=1} \frac{t_{d,s+d}t_{d,s+d+1}\cdots t_{d,r}}{t_{d+1,s+d}t_{d+1,s+d+1}\cdots t_{d+1,r}}}\nonumber\\
&=&a_{m+1}\cdots a_{m+d}\prod^d_{i=1}t_{i,i}t_{i,i+1}\cdots t_{i,s+i-1}\label{pr-3} \\
& & \frac{c (\sum^{m}_{s=1} t_{d+1,1}\cdots t_{d+1,s-1+d}
t_{d,s+d}t_{d,s+d+1}\cdots t_{d,r})
 + \sum^{r-d+1}_{s=m+1} t_{d+1,1}\cdots t_{d+1,s-1+d} t_{d,s+d}\cdots t_{d,r}}
{\sum^{r-d+1}_{s=1} t_{d+1,1}\cdots t_{d+1,s-1+d} t_{d,s+d}t_{d,s+d+1}\cdots t_{d,r}}.\nonumber
\end{eqnarray}
The relation (\ref{pr-1}) follows from (\ref{pr-2}) and (\ref{pr-3}). \qed

\subsection{An example}

\begin{ex}

We consider the case $SL_5(\mathbb{C})$, $\textbf{i}_0=(1,2,3,4,1,2,3,1,2,1)$.
The matrix $\tilde{B}_{\textbf{i}_0}$ is as follows:
\begin{eqnarray*}
\tilde{B}_{\textbf{i}_0}= \bordermatrix{
& -4 & -3 & -2 & -1 & 1 & 2 & 3 & 4 & 5 & 6 & 7 & 8 & 9 & 10 \cr
-4 & 1 & 0 & 0 & 0 & 0 & 0 & -1 & 1 & 0 & 0 & 0 & 0 & 0 & 0 \cr
-3 & -1 & 1 & 0 & 0 & 0 & -1 & 1 & 0 & 0 & 0 & 0 & 0 & 0 & 0 \cr
-2 & 0 & -1 & 1 & 0 & -1 & 1 & 0 & 0 & 0 & 0 & 0 & 0 & 0 & 0 \cr
-1 & 0 & 0 & -1 & 1 & 1 & 0 & 0 & 0 & 0 & 0 & 0 & 0 & 0 & 0 \cr
1 & 0 & 0 & 1 & -1 & 0 & -1 & 0 & 0 & 1 & 0 & 0 & 0 & 0 & 0 \cr
2 & 0 & 1 & -1 & 0 & 1 & 0 & -1 & 0 & -1 & 1 & 0 & 0 & 0 & 0 \cr
3 & 1 & -1 & 0 & 0 & 0 & 1 & 0 & -1 & 0 & -1 & 1 & 0 & 0 & 0 \cr
4 & -1 & 0 & 0 & 0 & 0 & 0 & 1 & 1 & 0 & 0 & -1 & 0 & 0 & 0 \cr
5 & 0 & 0 & 0 & 0 & -1 & 1 & 0 & 0 & 0 & -1 & 0 & 1 & 0 & 0 \cr
6 & 0 & 0 & 0 & 0 & 0 & -1 & 1 & 0 & 1 & 0 & -1 & -1 & 1 & 0 \cr
7 & 0 & 0 & 0 & 0 & 0 & 0 & -1 & 0 & 0 & 1 & 1 & 0 & -1 & 0 \cr
8 & 0 & 0 & 0 & 0 & 0 & 0 & 0 & 0 & -1 & 1 & 0 & 0 & -1 & 1 \cr
9 & 0 & 0 & 0 & 0 & 0 & 0 & 0 & 0 & 0 & -1 & 0 & 1 & 1 & -1 \cr
10 & 0 & 0 & 0 & 0 & 0 & 0 & 0 & 0 & 0 & 0 & 0 & -1 & 0 & 1 \cr
}
\end{eqnarray*}

Let $A:=(A_{-4},A_{-3},A_{-2},A_{-1},A_1,A_2,A_3,A_4,A_5,A_6,A_7,A_8,A_9,A_{10})\in \mathcal{A}_{\Sigma_{\textbf{i}_0}}=(\mathbb{C}^{\times})^{14}$, $M$ be the matrix in Theorem \ref{pgdef}, and $p=p_M:\mathcal{A}_{\Sigma_{\textbf{i}_0}}\rightarrow \mathcal{X}_{\Sigma_{\textbf{i}_0}}$ be the map in Proposition \ref{map-p}. We set $(P_{-4},P_{-3},P_{-2},P_{-1},P_1,\cdots,P_{10})=p(A)$. 

In this example, we shall calculate $(a_{\Sigma_{\textbf{i}_0}}^{-1}\circ (\iota\circ\zeta^{u,e})^{-1} \circ e_1^c \circ  (\iota\circ\zeta^{u,e}) \circ a_{\Sigma_{\textbf{i}_0}}A)_{k}$ via two ways: (1) An way using Theorem \ref{thm4}, (2) an way using Theorem \ref{thm2}.

First, let us calculate it using Theorem \ref{thm4}. Taking Theorem \ref{thm1} into account, we get
\begin{eqnarray*}
& &x_{\Sigma_{\bf{{\rm i}}}}^{-1}\circ e^c_1\circ x_{\Sigma_{\bf{{\rm i}}}}\circ p(A)\\
&=&\left(P_{-4},P_{-3},P_{-2}\frac{ P_1P_5P_8+P_5P_8+P_8+1}{cP_1P_5P_8+P_5P_8+P_8+1} ,
P_{-1}\frac{c P_1P_5P_8+P_5P_8+P_8+1}{P_1P_5P_8+P_5P_8+P_8+1},\right.\\
& &P_1\frac{c(P_1P_5P_8+P_5P_8)+P_8+1}{P_1P_5P_8+P_5P_8+P_8+1},
 P_2 \frac{c P_1P_5P_8+P_5P_8+P_8+1}{c(P_1P_5P_8+P_5P_8)+P_8+1},P_3,P_4,\\
& & P_5\frac{c(P_1P_5P_8+P_5P_8+P_8)+1}{cP_1P_5P_8+P_5P_8+P_8+1},P_6\frac{c (P_1P_5P_8+P_5P_8)+P_8+1}{c(P_1P_5P_8+P_5P_8+P_8)+1},\\
& & P_7,P_8\frac{c (P_1P_5P_8+P_5P_8+P_8+1)}{c(P_1P_5P_8+P_5P_8)+P_8+1},P_9\frac{c(P_1P_5P_8+P_5P_8)+P_8+1}{c P_1P_5P_8+P_5P_8+P_8+1},\\
& & \left. P_{10}\frac{c (P_1P_5P_8+P_5P_8+P_8+1)}{c(P_1P_5P_8+P_5P_8+P_8)+1}\right).
\end{eqnarray*}

Hence, we can calculate
\begin{eqnarray}
& &(a_{\Sigma_{\textbf{i}_0}}^{-1}\circ (\iota\circ\zeta^{u,e})^{-1} \circ e_1^c \circ  (\iota\circ\zeta^{u,e}) \circ a_{\Sigma_{\textbf{i}_0}}A)_{10^-} \nonumber \\
&=& (a_{\Sigma_{\textbf{i}_0}}^{-1}\circ (\iota\circ\zeta^{u,e})^{-1} \circ e_1^c \circ  (\iota\circ\zeta^{u,e}) \circ a_{\Sigma_{\textbf{i}_0}}A)_{8} \nonumber \\
&=& A_8 c^{D_{10,10}} \frac{P_{10}}{(x_{\Sigma_{\bf{{\rm i}}}}^{-1}\circ e^c_1\circ x_{\Sigma_{\bf{{\rm i}}}}\circ p(A))_{10}} \nonumber \\
&=& A_8 \left(\frac{c(P_1P_5P_8+P_5P_8+P_8)+1}{P_1P_5P_8+P_5P_8+P_8+1}\right).\label{ex-0}
\end{eqnarray}
Similarly, we obtain
\begin{eqnarray}
& &(a_{\Sigma_{\textbf{i}_0}}^{-1}\circ (\iota\circ\zeta^{u,e})^{-1} \circ e_1^c \circ  (\iota\circ\zeta^{u,e}) \circ a_{\Sigma_{\textbf{i}_0}}A)_{9^-}\nonumber\\
&=&(a_{\Sigma_{\textbf{i}_0}}^{-1}\circ (\iota\circ\zeta^{u,e})^{-1} \circ e_1^c \circ  (\iota\circ\zeta^{u,e}) \circ a_{\Sigma_{\textbf{i}_0}}A)_{6}\nonumber\\
&=& A_6 c^{D_{\{9,10\},\{8,10\}}}
\frac{P_{9}}{(x_{\Sigma_{\bf{{\rm i}}}}^{-1}\circ e^c_1\circ x_{\Sigma_{\bf{{\rm i}}}}\circ p(A))_{9}} 
\left(\frac{P_{10}}{(x_{\Sigma_{\bf{{\rm i}}}}^{-1}\circ e^c_1\circ x_{\Sigma_{\bf{{\rm i}}}}\circ p(A))_{10}}\right)^{D_{9,8}}\nonumber\\
&=&A_6, \label{ex-2}
\end{eqnarray}
and 
\begin{equation}\label{ex-5}
(a_{\Sigma_{\textbf{i}_0}}^{-1}\circ (\iota\circ\zeta^{u,e})^{-1} \circ e_1^c \circ  (\iota\circ\zeta^{u,e}) \circ a_{\Sigma_{\textbf{i}_0}}A)_5
=A_5 \left(\frac{c(P_1P_5P_8+P_5P_8)+P_8+1}{P_1P_5P_8+P_5P_8+P_8+1}\right),
\end{equation}
\begin{equation}\label{ex-6}
(a_{\Sigma_{\textbf{i}_0}}^{-1}\circ (\iota\circ\zeta^{u,e})^{-1} \circ e_1^c \circ  (\iota\circ\zeta^{u,e}) \circ a_{\Sigma_{\textbf{i}_0}}A)_3
=A_3,\ (a_{\Sigma_{\textbf{i}_0}}^{-1}\circ (\iota\circ\zeta^{u,e})^{-1} \circ e_1^c \circ  (\iota\circ\zeta^{u,e}) \circ a_{\Sigma_{\textbf{i}_0}}A)_2
=A_2,
\end{equation}
\begin{equation}\label{ex-7}
(a_{\Sigma_{\textbf{i}_0}}^{-1}\circ (\iota\circ\zeta^{u,e})^{-1} \circ e_1^c \circ  (\iota\circ\zeta^{u,e}) \circ a_{\Sigma_{\textbf{i}_0}}A)_1
=A_1 \left(\frac{cP_1P_5P_8+P_5P_8+P_8+1}{P_1P_5P_8+P_5P_8+P_8+1}\right),
\end{equation}
\begin{equation}\label{ex-8}
(a_{\Sigma_{\textbf{i}_0}}^{-1}\circ (\iota\circ\zeta^{u,e})^{-1} \circ e_1^c \circ  (\iota\circ\zeta^{u,e}) \circ a_{\Sigma_{\textbf{i}_0}}A)_{-j}
=A_{-j}\ \ (1\leq j\leq 4),
\end{equation}
\begin{equation}\label{ex-9}
(a_{\Sigma_{\textbf{i}_0}}^{-1}\circ (\iota\circ\zeta^{u,e})^{-1} \circ e_1^c \circ  (\iota\circ\zeta^{u,e}) \circ a_{\Sigma_{\textbf{i}_0}}A)_{k}
=c^{\delta_{k,10}}A_{k}\ \ (k=4,7,9,10).
\end{equation}


Next, by Theorem \ref{thm2}, we get
\begin{eqnarray*}
(a_{\Sigma_{\textbf{i}_0}}^{-1}\circ (\iota\circ\zeta^{u,e})^{-1} \circ e_1^c \circ  (\iota\circ\zeta^{u,e}) \circ a_{\Sigma_{\textbf{i}_0}}A)_{8} 
&=& (a_{\Sigma_{\textbf{i}_0}}^{-1}\circ (\iota\circ\zeta^{u,e})^{-1} \circ e_1^c \circ  (\iota\circ\zeta^{u,e}) \circ a_{\Sigma_{\textbf{i}_0}}A)_{3,1} \\
&=& A_{3,1} \frac{c(P_{1,1}P_{2,1}P_{3,1}+P_{2,1}P_{3,1}+P_{3,1})+1}{P_{1,1}P_{2,1}P_{3,1}+P_{2,1}P_{3,1}+P_{3,1}+1} \\
&=& A_{8} \frac{c(P_{1}P_{5}P_{8}+P_{5}P_{8}+P_{8})+1}{P_{1}P_{5}P_{8}+P_{5}P_{8}+P_{8}+1},
\end{eqnarray*}
\begin{eqnarray*}
(a_{\Sigma_{\textbf{i}_0}}^{-1}\circ (\iota\circ\zeta^{u,e})^{-1} \circ e_1^c \circ  (\iota\circ\zeta^{u,e}) \circ a_{\Sigma_{\textbf{i}_0}}A)_{6} 
&=& (a_{\Sigma_{\textbf{i}_0}}^{-1}\circ (\iota\circ\zeta^{u,e})^{-1} \circ e_1^c \circ  (\iota\circ\zeta^{u,e}) \circ a_{\Sigma_{\textbf{i}_0}}A)_{2,2} \\
&=& A_{2,2} =A_6,
\end{eqnarray*}
and these results coincide with (\ref{ex-0}) and (\ref{ex-2}).
Similarly, we can verify the results of calculations
for $(a_{\Sigma_{\textbf{i}_0}}^{-1}\circ (\iota\circ\zeta^{u,e})^{-1} \circ e_1^c \circ  (\iota\circ\zeta^{u,e}) \circ a_{\Sigma_{\textbf{i}_0}}A)_{k}$
by using Theorem \ref{thm2} coincide with (\ref{ex-5})-(\ref{ex-9}) ($k=-4,-3,-2,-1,1,2,3,4,5,7,9,10$).

\end{ex}



\end{document}